\documentclass[a4paer,10pt]{article}
\usepackage{soul}
\usepackage{pifont}
\usepackage{rotating}
\usepackage[doc]{optional}
\usepackage{xcolor}
\usepackage{exscale,relsize}
\usepackage{amsmath}
\usepackage{amsfonts}
\usepackage{amssymb}
\usepackage[mathscr]{euscript}
\usepackage{mathrsfs}
\usepackage{calc}
\usepackage{theorem}
\usepackage{pifont}
\usepackage{graphicx}
\definecolor{labelkey}{rgb}{0,0.08,0.45}
\definecolor{rekey}{rgb}{0,0.6,0.0}
\definecolor{Brown}{rgb}{0.45,0.0,0.05}
\newtheorem{theorem}{Theorem}[section]
\newtheorem{lemma}[theorem]{Lemma}

\newtheorem{corollary}[theorem]{Corollary}

\newtheorem{proposition}[theorem]{Proposition}
\newtheorem{definition}[theorem]{Definition}

\theoremstyle{plain}{\theorembodyfont{\rmfamily}
}
\theoremstyle{plain}{\theorembodyfont{\rmfamily}
}
\theoremstyle{plain}{\theorembodyfont{\rmfamily}
}
\theoremstyle{plain}{\theorembodyfont{\rmfamily}
\newtheorem{example}[theorem]{Example}}
\theoremstyle{plain}{\theorembodyfont{\rmfamily}
\newtheorem{remark}[theorem]{Remark}}

\theoremstyle{plain}{\theorembodyfont{\rmfamily}
}
\def\disp{\displaystyle}

\def\hat{\widehat}
\def\Tilde{\widetilde}

\def\Bar{\overline}
\def\ra{\rangle}
\def\la{\langle}
\def\ve{\varepsilon}
\def\B{\mathbb{B}}

\def\R{\mathbb{R}}

\def\ox{\bar{x}}

\def\oy{\bar{y}}

\def\ou{\bar{u}}

\def\kk{\kappa}

\def\bd{\mbox{\rm bd}\,}

\def\bd{\mbox{\rm bd}\,}

\def\cone{\mbox{\rm cone}\,}

\def\O{\Omega}

\def\emp{\emptyset}

\def\lm{\lambda}
\def\gg{\gamma}
\def\dd{\delta}
\def\al{\alpha}

\newcommand{\NN}{\ensuremath{\mathbb N}}

\newcommand{\di}{\ensuremath{\operatorname{dist}}}

\newcommand{\qed}{\hspace*{\fill}$\Box$\medskip}

\begin{document}

\title{\textsf{Error Bounds for Parametric Polynomial Systems with Applications to Higher-Order Stability Analysis and Convergence Rates}}
\author{G. LI\thanks{Department of Applied Mathematics, University of New South Wales, Sydney 2052, Australia. E-mail: g.li@unsw.edu.au. Research of this author was partly supported by the Australian Research Council.}, \quad B. S. MORDUKHOVICH\thanks{Department of Mathematics, Wayne State University, Detroit, MI 48201, USA. Email: boris@math.wayne.edu. Research of this author was partly supported by the USA National Science Foundation under grants DMS-1007132 and DMS-1512846 and by the USA Air Force Office of Scientific Research under grant No.\,15RT0462.}, \quad T. T. A. NGHIA\thanks{Department of Mathematics and Statistics, Oakland University, Rochester, MI 48306, USA. E-mail: nttran@oakland.edu.}, \quad and \quad T. S. PH\d{A}M\thanks{Department of Mathematics, University of Dalat, Dalat, Vietnam. E-mail: sonpt@dlu.edu.vn. Research of this author was partly supported by the Vietnam National Foundation for Science and Technology Development (NAFOSTED) under grant No.\,101.04-2013.07.}}
\maketitle
\vspace{-1.1cm}
\begin{center}
{\bf Dedicated to Terry Rockafellar in honor of his 80th birthday}
\end{center}
\begin{abstract} The paper addresses parametric inequality systems described by polynomial functions in finite dimensions, where state-dependent infinite parameter sets are given by finitely many polynomial inequalities and equalities. Such systems can be viewed, in particular, as solution sets to problems of generalized semi-infinite programming with polynomial data. Exploiting the imposed polynomial structure together with powerful tools of variational analysis and semialgebraic geometry, we establish a far-going extension of the \L ojasiewicz gradient inequality to the general nonsmooth class of supremum marginal functions as well as higher-order (H\"older type) local error bounds results with explicitly calculated exponents. The obtained results are applied to higher-order quantitative stability analysis for various classes of optimization problems including generalized semi-infinite programming with polynomial data, optimization of real polynomials under polynomial matrix inequality constraints, and polynomial second-order cone programming. Other applications provide explicit convergence rate estimates for the cyclic projection algorithm to find common points of convex sets described by matrix polynomial inequalities and for the asymptotic convergence of trajectories of subgradient dynamical systems in semialgebraic settings.


\end{abstract}
\vspace{-0.6cm}
\vspace{-0.3cm}

\section{Introduction}\vspace*{-0.1in}
\setcounter{equation}{0}

This paper is largely devoted to {\em polynomial semi-infinite optimization} and related topics being revolved around deriving {\em explicit error bounds} for infinite parametric inequality systems with (real) {\em polynomial data} as well as their various applications to stability analysis in optimization and convergence of algorithms. The imposed polynomial structure allows us to widely use powerful tools of semialgebraic geometry, while parametric inequalities naturally call for applying constructions and results of variational analysis and generalized differentiation. Needless to say that the seminal contributions by Terry Rockafellar to variational analysis and optimization are difficult to overstate, and it is our honor to dedicate this paper to him.\vspace*{-0.1in}

The primary attention of this paper is paid to the {\em parametric inequality systems}
\begin{equation}\label{eq:S}
S:=\big\{x\in\mathbb{R}^n\big|\;f_l(x,y)\le 0\;\mbox{ for all }\;y\in Y(x),\;l=1,\ldots,L\big\},
\end{equation}
where each function $f_l\colon\mathbb{R}^n\times\mathbb{R}^m$ as $l\in L$ for the given natural number $L\in\mathbb{N}$ is a {\em polynomial}, and where $Y:\R^n\rightrightarrows\R^m$ is a set-valued mapping that is also described by finitely many polynomials via inequality and equality constraints. Systems of type \eqref{eq:S} naturally arise as feasible solution sets in problems of generalized semi-infinite programming, second-order cone programming, robust optimization, and matrix inequalities with polynomial data; see below for more details and applications.\vspace*{-0.1in}

One of the most important issues associated with the inequality systems (\ref{eq:S}) is establishing the so-called {\em error bounds}. Given $\bar x\in\mathbb{R}^n$, recall that a (local) {\em error bound} of $S$ with a {\em H\"older exponent} $\tau\in(0,1]$ holds for $F$ at $\ox$ if there exist
constants $c>0$ and $\ve>0$ such that
\begin{equation}\label{eq:error_bound}
{\rm dist}(x,S)\le c\,\bigg(\Big[\sup_{y\in Y(x),\atop 1\le l\le L}f_l(x,y)\Big]_{+}\bigg)^{\tau}\;\mbox{ for all }\;x\;\mbox{ with }\;\|x-\overline{x}\|\le\ve,
\end{equation}
where ${\rm dist}(x,S)$ signifies the Euclidean distance between $x$ and $S$, and where $[\alpha]_{+}:=\max\{\alpha,0\}$. The supremum in \eqref{eq:error_bound} is obviously achieved and it can be replaced by `max' if $Y(x)$ is closed and bounded.\vspace*{-0.1in}

The study of error bounds has attracted a lot of attention of many researchers over the years and has found numerous applications to, in particular, sensitivity analysis for various problems of mathematical programming, termination criteria for numerical algorithms, etc. We refer the reader to \cite{Pang1997} for an excellent survey in these directions and to the more recent papers \cite{fhko,Klatte,kr,Boris,Li1,Thera,Ye} with the bibliographies therein. It is worth noting that the major attention in the aforementioned and many other publications on error bounds has drawn to the case of {\em linear rate} ($\tau=1)$, where this issue is related to metric regularity and subregularity notions in basic variational analysis. Our main interest in this paper concerns {\em fractional/root} error bounds $0<\tau<1$ in \eqref{eq:error_bound}. For the case of {\em finite} and {\em fixed} sets $Y(x)\equiv\Omega$ in  \eqref{eq:error_bound}, some results in this direction have been obtained in, e.g., \cite{Borwein,Dinh2014-1,Li_MP,Boris,Luo_Luo,Luo,Ng} with various applications therein.\vspace*{-0.1in}

It is proved \cite{Luo_Luo} in this finite case of $Y(x)\equiv\O$, by using the cerebrated \L ojasiewicz gradient inequality \cite{lo}, that (\ref{eq:error_bound}) holds with some {\em unknown} exponent $\tau\in(0,1)$ for polynomial systems \eqref{eq:S}. Employing advanced techniques of variational analysis, we have recently derived in this case \cite{LMP2014} several error bounds with exponents {\em explicitly determined} by the dimension of the underlying space and the number/degree of the involved polynomials. The techniques and results developed in \cite{LMP2014} allowed us to resolve several open questions raised in the literature, which include establishing explicit H\"{o}lder error bounds for nonconvex quadratic systems and higher-order semismoothness of the maximum eigenvalue for symmetric tensors.\vspace*{-0.1in}

The primary goal of this paper is to obtain explicit error bounds of type \eqref{eq:error_bound} for polynomial inequality systems with {\em infinite} and {\em variable} sets $Y(x)\subset\R^m$. Besides undoubted importance of these issues for their own sake, we have been motivated by applications of infinite polynomial systems and error bounds for them to higher-order stability and convergence rates of algorithms in optimization-related areas as well as in asymptotic analysis of dynamical systems, where estimates of type \eqref{eq:error_bound} with infinite sets $Y(x)$ are crucial.\vspace*{-0.1in}

As the reader can see below, deriving error bounds for the case of infinite and variable sets $Y(x)$ in \eqref{eq:S} is significantly more involved in comparison with our developments for finite systems in \cite{LMP2014}. First we present the following three-dimensional example showing that the error bound \eqref{eq:error_bound} may fail for any $\tau>0$ for infinite polynomial inequality systems even in the case of constant sets $Y(x)\equiv\O$.\vspace*{-0.1in}

\begin{example}{\bf(failure of H\"older error bounds for infinite polynomial systems).}\label{ex:1} Consider the polynomial system of type \eqref{eq:S} containing only one inequality given by $f:\mathbb{R}\times\mathbb{R}^2\rightarrow\mathbb{R}$ in the form
\begin{equation*}
f(x,y)\le 0\;\mbox{ with }\;y\in\O,\;\mbox{ where }\;f(x,y):=x y_1+y_2,\;y=(y_1,y_2),
\end{equation*}
where the infinite set $\O$ is constructed as follows. Take the ${\cal C}^\infty$-smooth function of one variable
\begin{equation*}
\phi(x):=\left\{\begin{array}{ccc}
e^{-\frac{1}{x^2}}&\mbox{ if }&x\ne 0,\\
0&\mbox{ if }&x=0
\end{array}\right.
\end{equation*}
and define the set $\O\subset\R^2$ by the conditions
\begin{equation}\label{exO}
\Omega:=\big\{(y_1,y_2)\big|\;\exists\,t\in[-0.5,0.5]\;\mbox{ such that }\;y_1=\phi'(t)\;\mbox{ and }\;y_2=\phi(t)-t\phi'(t)\big\},
\end{equation}
where $\phi'$ stands for the derivative of $\phi$. Since $\phi$ is ${\cal C}^{\infty}$-smooth, the set $\Omega$ in \eqref{exO} is nonempty and compact as the image of a compact interval under a continuous mapping. We claim that
\begin{equation}\label{eq:claim1}
\max_{y\in\Omega}f(x,y)=\phi(x)\;\mbox{ whenever }\;x\in[-0.5,0.5].
\end{equation}
Indeed, for any $x\in[-0.5,0.5]$ and $y\in\Omega$ there is $t\in[-0.5,0.5]$ with $y_1=\phi'(t)$ and $y_2=\phi(t)-t\phi'(t)$. Thus
\[
f(x,y)=x\phi'(t)+\phi(t)-t\phi'(t)=\phi(t)+(x-t)\phi'(t).
\]
Applying the second-order Taylor expansion to the function $\phi$ tells us that
\[
\phi(x)=\phi(t)+\phi'(t)(x-t)+\frac{\phi''(a)}{2}(x-t)^2,\;\mbox{ for some }\;a\in[x,t]\subset[-0.5,0.5].
\]
Note that $\phi''(a)=(4a^{-6}-6a^{-4})e^{-\frac{1}{a^2}}=2a^{-6}e^{-\frac{1}{a^2}}\left(2-3a^2\right)\ge 0$, where the last inequality holds due to the choice of $a\in[-0.5,0.5]$. Hence we get the relationships
\[
f(x,y)=\phi(t)+(x-t)\phi'(t)\le\phi(x)\;\mbox{ for all }\;x\in[-0.5,0.5],\;y\in\Omega,
\]
which imply the inequality $\max_{y\in\Omega}f(x,y)\le\phi(x)$ whenever $x\in[-0.5,0.5]$. On the other hand, it follows from the above constructions of $f$ and $\Omega$ that $$
\max_{y\in\Omega}f(x,y)\ge x\phi'(x)+\big(\phi(x)-x\phi'(x)\big)=\phi(x)\;\mbox{ for each }\;x\in[-0.5,0.5],
$$
which therefore justifies the claim in \eqref{eq:claim1}. Having this in mind, consider the set
\begin{equation}\label{eq:S1}
S:=\big\{x\in\mathbb{R}\big|\;f(x,y)\le 0\;\mbox{ for all }\;y\in\Omega\big\}
\end{equation}
and observe that $S\cap[-0.5,0.5]=\{x|\;\phi(x)\le 0\}=\{0\}$. Let us  now check that the local H\"{o}lder error bound \eqref{eq:error_bound} fails for \eqref{eq:S1} at $\ox=0$ for any exponent $\tau>0$. To see it, take $x_k:=\frac{1}{\sqrt{\ln(k+1)}}\rightarrow 0$ as $k\in\mathbb{N}$ and get for all large $k$ that
$d\left(x_k,S\right)=d\left(x_k,S\cap[-0.5,0.5]\right)=|x_k|$. This allows us to conclude that
\[
\frac{d(x_k,S)}{[\max_{y\in\Omega}f(x_k,y)]_+^{\frac{1}{\tau}}}=\frac{|x_k|}{[\phi(x_k)]_+^{\frac{1}{\tau}}}=\frac{(k+1)^{\frac{1}{\tau}}}{\sqrt{\ln (k+1)}} \rightarrow\infty\;\mbox{ as }\;k\rightarrow\infty
\]
whenever $\tau>0$ is chosen, and thus the error bound \eqref{eq:error_bound} fails for the system $S$ in \eqref{eq:S1}.
\end{example}\vspace*{-0.2in}

In what follows we prove that such a situation {\em does not emerge} if the sets $Y(x)$ in \eqref{eq:S} are described by
\begin{equation}\label{Y-poly}
Y(x)=\big\{y\in\mathbb{R}^m\big|\;g_i(x,y)\le 0,\;i=1,\ldots,r,\;\mbox{ and }\;h_j(x,y)=0,\;j=1,\ldots,s\big\},
\end{equation}
where $g_i$ and $h_j$ are {\em polynomials}. It is shown in Section~4 that $\tau$ in \eqref{eq:error_bound} is {\em explicitly calculated} in terms of degrees of the polynomials and dimensions of the spaces in question. The key of our analysis is a new nonsmooth extension of the {\em \L ojasiewicz inequality} to the class of {\em supremum marginal functions}
\begin{equation*}
\phi(x):=\sup_{y\in Y(x)}f(x,y)
\end{equation*}
described by polynomials $f$ and $g_i,h_j$ in \eqref{Y-poly}. This is done in Section~3 by using powerful tools of variational analysis and semialgebraic geometry reviewed in Section~2.\vspace*{-0.1in}

Sections~5 and 6 are devoted to applications. In Section~5 we develop {\em quantitative higher-order stability analysis} for remarkable classes of polynomial optimization problems: {\em generalized semi-infinite programming}, {\em optimization of matrix inequalities}, and {\em second-order cone programming}. Finally, Section~6 contains explicit estimates of {\em convergence rates} for the {\em cyclic projection algorithm} to solve feasibility problems for convex sets described by matrix polynomial inequalities and also for asymptotic analysis of {\em subgradient dynamical systems} governed by maximum functions with polynomial data.\vspace*{-0.2in}

\section{Tools of Variational Analysis and Semialgebraic Geometry}\label{SectionPreliminary}\vspace*{-0.1in}
\setcounter{equation}{0}

This section briefly discusses some tools of generalized differentiation in variational analysis and of semialgebraic geometry widely used in the paper. Throughout this work we deal with finite-dimensional Euclidean spaces labeled as $\mathbb{R}^n$ and endowed by the inner product $\langle x,y\rangle=x^Ty$. The symbol $\mathbb{B}_\ve(x)$ (resp.\ $\overline{\mathbb{B}}_\ve(x)$) stands for the open (resp.\ closed) ball with center $x$ and radius $\ve>0$ while $\mathbb{B}$ (resp.\ $\overline{\mathbb{B}}$) stands for the open (resp closed) unit ball centered at the origin in $\mathbb{R}^n$. Given a set $A\subset\mathbb{R}^n$, its interior (resp.\  boundary, convex hull, and conic convex hull) is denoted by ${\rm int}\,A$ (resp.\ ${\rm bd}\,A$, ${\rm co}\,A$, and ${\rm cone}\,A$).\vspace*{-0.1in}

Starting with {\em variational analysis}, recall first two subdifferential notions needed in what follows. The reader can find more information and references in the books \cite{Boris_Book,Rockafellar1998}.\vspace*{-0.1in}

Given a function $f\colon{\Bbb R}^n\rightarrow{\Bbb R}$ continuous around $x$, the {\em proximal subdifferential} of $f$ at $x$ is
\begin{equation}\label{PS}
\partial_P f(x):=\left\{v\in{\Bbb R}^n\Big|\;\liminf_{\|h\|\to 0,\ h\ne 0}\frac{f(x+h)-f(x)-\langle v,h\rangle}{\|h\|^2}>-\infty\right\}.
\end{equation}
The {\em limiting subdifferential} of $f$ at $x$ (known also as the general, basic or Mordukhovich subdifferential) is
\begin{equation}\label{ls}
\partial f(x):=\big\{v\in\R^n\big|\;\exists\,x_k\to x,\;v_k\to v\;\mbox{ with }\;v_k\in{\partial}_P f(x_k),\;k\in\mathbb{N}\big\}.
\end{equation}
We clearly have $\partial_P f(x)\subset\partial f(x)$, where the first set may often be empty (while not so in a dense sense), but the second one is nonempty for any locally Lipschitzian function. Furthermore, the set $\partial_P f(x)$ is always convex but may not be closed, while $\partial f(x)$ is closed but may often be nonconvex. Both subdifferentials \eqref{PS} and \eqref{ls} reduce to the gradient $\{\nabla f(x)\}$ for smooth functions and to the subdifferential of convex analysis for convex ones. A significant advantage of the limiting subdifferential \eqref{ls} is {\em full calculus} in the general nonconvex setting that is based on variational and extremal principles; see \cite{Boris_Book,Rockafellar1998} for more details.\vspace*{-0.1in}

The major variational notion used below is the {\em limiting subdifferential slope} of $f$ at $x$ defined via \eqref{ls} by
\begin{equation}\label{slope}
\frak m_f(x):=\inf\big\{\|v\|\big|\;v\in\partial f(x)\big\},
\end{equation}
where $\inf\{\emp\}:=\infty$. It reduces to the classical gradient slope $\frak m_f(x)=\|\nabla f(x)\|$ for smooth functions.\vspace*{-0.08in}

Let us next formulate some continuity notions for set-valued mappings $F:\mathbb{R}^n\rightrightarrows\mathbb{R}^m$; see, e.g., \cite{Rockafellar1998}. It is said that $F$ is {\em outer semicontinuous} (o.s.c.) at $\ox$ if for any sequence $(x_k,y_k)$ converging to $(\ox,\oy)$ with $y_k\in F(x_k)$ we have $\oy\in F(\ox)$.
Further, $F$ is {\em inner semicontinuous} (i.s.c.) at $\ox$ if for any sequence $x_k\to\ox$ and any $\oy\in F(\ox)$ there are $y_k\in F(x_k)$ as $k\in\mathbb{N}$ satisfying $y_k\to\oy$. We say that $F$ is o.s.c.\ or i.s.c.\ {\em around} $\ox$ if it has this property at every $x$ in a neighborhood of $\ox$.
\vspace*{-0.1in}

Finally, we present some notions and facts from (real) {\em semialgebraic geometry} following \cite{Bochnak1998}. It is said that:\vspace*{-0.1in}

$\bullet$ $A\subset\mathbb{R}^n$ is a {\em semialgebraic set} if it is a finite union of subsets given by
$$
\big\{x\in\mathbb{R}^n\big|\;f_i(x)=0,\;i=1,\ldots,q,\;\mbox{ and }\;f_i(x)>0,\;i=q+1,\ldots,p\big\},
$$
where all the functions $f_i$, $i=1,\ldots,p$, are polynomials of some degrees.\vspace*{-0.1in}

$\bullet$ $f\colon A\to B$ is a {\em semialgebraic mapping} if it maps one semialgebraic set $A\subset\mathbb{R}^n$ to another one $B\subset\mathbb{R}^m$ and its graph $\{(x,y)\in A\times B|\;y=f(x)\}$ is a semialgebraic subset of $\Bbb{R}^n\times\Bbb{R}^m$. We say that $f\colon A\to B$ is {\em locally semialgebraic} around $\ox$ if there exists a neighborhood $V$ of the point $(\ox,f(\ox))$ such that the set $V\cap\{(x,y)\in A\times B|\;y=f(x)\}$ is semialgebraic in $\Bbb{R}^n\times\Bbb{R}^m$.\vspace*{-0.1in}

The class of semialgebraic sets is closed under taking finite intersections, finite unions, and complements; furthermore, a Cartesian product of semialgebraic sets is semialgebraic. A major fact concerning the class of semialgebraic sets is given by the following seminal result of semialgebraic geometry.\\[1ex]
{\bf Tarski-Seidenberg Theorem}. {\em Images of semialgebraic sets under semialgebraic maps are semialgebraic}.\vspace*{-0.1in}

We also need another fundamental result taken from \cite[Theorem~4.2]{Acunto2005}, which provides an exponent estimate in the classical \L ojasiewicz gradient inequality for polynomials. For brevity we label it as:\\[1ex]
{\bf \L ojasiewicz Gradient Inequality}. {\em Let $f$ be a polynomial on $\mathbb{R}^n$ with degree $d\in\mathbb{N}$. Suppose that $f(0)=0$ and $\nabla f(0)=0$. Then there exist constants $c,\ve>0$ such that for all $x\in\mathbb{R}^n$ with $\|x\|\le\ve$ we have
\begin{equation}\label{loj}
\|\nabla f(x)\|\ge c|f(x)|^{1-\tau},\;\mbox{ where }\;\tau ={\mathscr R}(n,d)^{-1}\;\mbox{ and }\;{\mathscr R}(n,d):=\begin{cases}
1&\mbox{ if }\;d=1,\\
d(3d-3)^{n-1}&\mbox{ if }\;d\ge 2.
\end{cases}
\end{equation}}
\vspace*{-0.4in}

\section{\L ojasiewicz Inequality for Supremum Marginal Functions}\vspace*{-0.1in}
\setcounter{equation}{0}

The main aim of this section is extending the \L ojasiewicz gradient inequality \eqref{loj} to the following class of (polynomial) {\em supremum marginal functions} given by
\begin{equation}\label{max}
\phi(x):=\disp\sup_{y\in Y(x)}f(x,y),
\end{equation}
where $f\colon\R^n\times\R^m\to\R$ is a polynomial of degree at most $d$, and where the set-valued mapping $Y\colon\R^n\rightrightarrows\R^m$ is defined in \eqref{Y-poly} by polynomials $g_i\colon\R^n\times\R^m\to\R$ for $i=1,\ldots,r$ and $h_j\colon\R^n\times\R^m\to\R$ for $j=1,\ldots,s$ of degree at most $d$. Functions of type \eqref{max} are intrinsically nonsmooth, and thus deriving a nonsmooth counterpart of \eqref{loj} for them requires the usage of an appropriate subdifferential of $\phi$. The reader will see below that a nonsmooth version of the \L ojasiewicz inequality for $\phi$ in terms of the limiting subdifferential slope $\frak m_\phi(x)$ from \eqref{slope}, which replaces the gradient norm in \eqref{loj}, plays a key role in establishing our H\"{o}lder-type local error bounds and their subsequent applications in this paper.\vspace*{-0.1in}

To proceed in this direction, we have to calculate the limiting subdifferential of $\phi$ from \eqref{max} in terms of its initial data, which is not an easy task by taking into account that the sets $Y(x)$ are infinite and variable. When $Y(x)\equiv\O$, \eqref{max} reduces to the {\em supremum function} for which the most recent subdifferential results can be found in \cite{MN13}; see also the references therein. When the supremum in \eqref{max} is replaced by the infimum, we arrive at the class of (infimum) {\em marginal functions} well investigated in variational analysis \cite{Boris_Book}. Needless to say that the supremum and infimum operations are essentially different in unilateral analysis and that (lower) subdifferential properties under supremum are significantly more involved and challenging.\vspace*{-0.1in}

We first show that the functions of type \eqref{max} enjoy the following important properties.\vspace*{-0.1in}
\begin{proposition}{\bf (supremum marginal functions are locally semialgebraic and upper semicontinuous).}\label{Lemma31} Let the sets $Y(x)$ in \eqref{Y-poly} be nonempty and uniformly bounded around some point $\ox\in\R^n$. Then the function $\phi$ in \eqref{max} is locally semialgebraic and locally upper semicontinuous $($u.s.c.$)$ around $\ox$.
\end{proposition}\vspace*{-0.15in}
\begin{proof}
It follows from the Tarski-Seidenberg Theorem and the assumptions made that $\phi$ is well defined and semialgebraic around $\ox$. Let us verify that $\phi$ is u.s.c.\ around $\ox$. Select a number $\ve>0$ and a compact set $C\subset\mathbb{R}^m$ so that $\emp\ne Y(x)\subset C$ for all $x\in\overline\B_\ve(\ox)$. We claim that $\phi$ is locally u.s.c.\ on $\B_\ve(\ox)$. Indeed, assume on the contrary that $\phi$ is not u.s.c.\ at some $x\in\B_\ve(\ox)$ and then get a sequence $x_k\to x$ and a number $\delta>0$ such that $\phi(x_k)>\phi(x)+\delta$. Using the continuity of $f$, the closedness of $Y(x)$ for all $x$, and the inclusion $Y(x_k)\subset C$ allows us to find $y_k\in Y(x_k)$ with $\phi(x_k)=f(x_k,y_k)$. By passing to a subsequence we may suppose without loss of generality that $y_k\to y$. It is easy to check that $y\in Y(x)$ due to the continuity of $g_i$ as $i=1,\ldots,r$ and $h_j$ as $j=1,\ldots,s$ in \eqref{Y-poly}. Since $f$ is continuous, we have $\phi(x_k)=f(x_k,y_k)\to f(x,y)\le\phi(x)$, which contradicts the assumption that $\phi(x_k)>\phi(x)+\delta$. Thus $\phi$ is locally u.s.c.\ around $\ox$.\qed
\end{proof}\vspace*{-0.15in}

Next we study the lower semicontinuity (l.s.c.) and H\"older continuity of the function $\phi$ \eqref{max} under rather mild assumptions. Define the {\em argmaximum set} important for our further analysis by
\begin{equation}\label{Y0}
M(x):=\big\{y\in Y(x)\big|\;\phi(x)=f(x,y)\big\}
\end{equation}
and recall that $\di(A;B):=\inf\{\|a-b\|\big|\;a\in A,\;b\in B\}$ for $A,B\subset\mathbb{R}^m$.
\vspace*{-0.1in}
\begin{proposition}{\bf (continuity of supremum marginal functions).}\label{Con} Given $\ox\in\mathbb{R}^n$, suppose that the sets $Y(x)$ are nonempty and  uniformly bounded around $\ox$. The following assertions hold:\\
{\bf (i)} If $\di(M(\ox);Y(x))\to 0$ as $x\to\ox$, then $\phi$ is l.s.c.\ $($and hence continuous$)$ at $\ox$. Thus the i.s.c.\ property of $Y$ at $\ox$ ensures that $\phi$ is continuous at this point.\\
{\bf (ii)} If there are numbers $\ve,L>0$ and $\al\in(0,1]$ such that
\begin{equation}\label{eq:00}
\di\big(M(x_1);Y(x_2)\big)\le L\|x_1-x_2\|^\al\;\mbox{ for all }\;x_1,x_2\in\B_\ve(\ox),
\end{equation}
then $\phi$ is H\"older continuous on $\B_\ve(\ox)$ with order $\al$.
\end{proposition}\vspace*{-0.2in}
\begin{proof}
Since the sets $Y(x)$ are nonempty and uniformly bounded around $\ox$, there exist a compact set $C$ and a number $\ve>0$ such that $Y(x)\subset C$ for all $x\in \B_\ve(\ox)$. Let us first verify (i) while arguing by contradiction. Suppose that $\di(M(\ox);Y(x))\to 0$ as $x\to\ox$ and that $\phi$ is not l.s.c.\ at $\ox$. Then there exist a constant $\nu>0$ and a sequence $x_k\to\ox$ such that $\phi(x_k)<\phi(\ox)-\nu$. Since $d(M(\ox);Y(x_k))\to 0$ as $k\rightarrow\infty$, we find $y_k\in Y(x_k)$ and $\oy_k\in M(\ox)$ with $\|y_k-\oy_k\|\to 0$. Using the boundedness of $\{\oy_k\}_{k\in \mathbb{N}}$ and passing to a subsequence if necessary allow us to find $\oy\in\R^m$ with $\oy_k\to\oy$ and thus $y_k\to\oy$. It is easy to see that $\oy\in Y(\ox)$ and that $\phi(\ox)=f(\ox,\oy_k)\to f(\ox,\oy)$, which yield in turn that $\oy\in M(\ox)$. Hence we have $\phi(x_k)\ge f(x_k,y_k)\to f(\ox,\oy)=\phi(\ox)$ as $k\to\infty$ while contradicting the assumption on $\phi(x_k)<\phi(\ox)-\nu$. This verifies that $\phi$ is l.s.c.\ at $\ox$ and hence its continuity there by Proposition~\ref{Lemma31}. The second statement in (i) follows from the first one and the definitions.\vspace*{-0.1in}

To justify now assertion (ii), suppose that \eqref{eq:00} holds and pick any $x_1,x_2\in \B_\ve(\ox)$. Since the sets $M(x)$ and $Y(x)$ are compact for all
$x\in\B_\ve(\ox)$, we deduce from \eqref{eq:00} that there are $y_1\in M(x_1)$ and $y_2\in Y(x_2)$ satisfying $\|y_1-y_2\|\le L\|x_1-x_2\|^\al$. Taking into account the Lipschitz continuity of the polynomial $f$ on $\Bar\B_\ve(\ox)\times C$ with some constant $\ell>0$, it follows that
\[\begin{array}{ll}
\phi(x_1)-\phi(x_2)&\disp\le f(x_1,y_1)-f(x_2,y_2)\le\ell(\|x_1-x_2\|+\|y_1-y_2\|)\\
&\disp\le\ell\big(\|x_1-x_2\|+L\|x_1-x_2\|^\al\big)=\ell\big((2\ve)^{1-\al}+L\big)\|x_1-x_2\|^\al.
\end{array}
\]
Similarly we conclude that $\phi(x_2)-\phi(x_1)\le\ell\big((2\ve)^{1-\al}+L\big)\|x_1-x_2\|^\al$, which verifies the H\"older continuity of $\phi$ on $\B_\ve(\ox)$ with order $\al$ and thus completes the proof of the proposition.\qed
\end{proof}\vspace*{-0.15in}
\begin{remark}{\bf (effective conditions for continuity of marginal functions).}\label{cont-mar} {\rm Proposition~\ref{Con}(i) tells us that $\phi$ in \eqref{max} is continuous at $\ox$ provided that $Y$ is i.s.c.\ at this point, which surely holds if $Y$ is {\em locally Lipschitzian} around $\ox$ in the standard/Hausdorff sense with taking into account that the sets $Y(x)$ from \eqref{Y-poly} are closed, nonempty, and uniformly bounded. Necessary and sufficient conditions for the local Lipschitz property of $Y$ around $\ox$ are given in \cite[Corollary~4.39]{Boris_Book} in terms of the gradients of the constraint functions $g_i,h_j$ at $(\ox,\oy)$ for {\em any} $\oy\in Y(\ox)$. In this way we get effective conditions ensuring the continuity of $\phi$ at (in fact {\em around}) $\ox$ via the initial data of \eqref{Y-poly}. On the other hand, \cite[Theorem~3.38(iv)]{Boris_Book} justifies the {\em Lipschitz continuity} of $\phi$ around $\ox$ under the i.s.c.\ property of the argmaximum mapping \eqref{Y0} around $(\ox,\oy)$ with some {\em fixed} $\oy\in M(\ox)$ and the {\em Lipschitz-like} (Aubin, pseudo-Lipschitz) property of $Y$ around this pair $(\ox,\oy)$, which is effectively characterized in \cite[Corollary~4.39]{Boris_Book} via the initial data.}
\end{remark}\vspace*{-0.2in}

To proceed further, we need the following qualification condition imposed at a reference point relative to some (in fact optimal) subset of the constraint set $Y(x)$ in \eqref{max}.\vspace*{-0.1in}

\begin{definition}{\bf(marginal constrained qualification).}\label{mmfcq} Given $\ox\in\R^n$ and $Y(\ox)$ from \eqref{Y-poly}, we say that the {\sc marginal Mangasarian-Fromovitz constraint qualification (MMFCQ)} holds at $\ox$ relative to some subset $\O\subset Y(\ox)$ if there is a vector $\xi\in \R^n$ such that
\begin{equation}\label{MF}
\left[\sum_{i=1}^r\lambda_i\nabla_xg_i(\ox,y)+\sum_{j=1}^s\kappa_j\nabla_x h_j(\ox,y)\right]^T\xi>0\;\mbox{ for any }\;y\in\O
\end{equation}
whenever Lagrange multipliers $(\lm,\kappa)\in\R_+^r\times\R^s\setminus\{0\}$ satisfy the conditions
\begin{equation}\label{CQ}
\sum_{i=1}^r\lambda_i\nabla_yg_i(\ox,y)+\sum_{j=1}^s\kappa_j\nabla_y h_j(\ox,y)=0\;\mbox{ and }\;\lm_i g_i(\ox,y)=0,\quad i=1,\ldots,r.
\end{equation}
\end{definition}\vspace*{-0.2in}

\begin{remark}{\bf (relationships of MMFCQ with EMFCQ and MFCQ).}\label{emfcq} {\rm The defined MMFCQ for the marginal supremum functions \eqref{max} is motivated by the {\em extended Mangasarian-Fromovitz constraint qualification} (EMFCQ) introduced in \cite{Jongen} for {\em generalized semi-infinite programs} (GSIPs) discussed below in Section~5. Reformulating EMFCQ for \eqref{max}, we see that it involves the cost function $f$ and also requires the validity of \eqref{MF} for a larger set of Lagrange multipliers in comparison with MMFCQ; so MMFCQ is weaker than EMFCQ in general. Observe also that when the functions $g_i$ and $h_j$ do not depend on $x$, the formulated MMFCQ condition means that there is no $(\lm,\kk)\in\R^r_+\times\R^s\setminus\{0\}$ satisfying \eqref{CQ}, which is equivalent to the conventional {\em Mangasarian-Fromovitz constraint qualification} (MFCQ) on $\O$ in the sense that
\begin{equation}\label{MFCQ}
\left.\begin{array}{r}\displaystyle\sum_{i=1}^r\lambda_i\nabla g_i(y)+\sum_{j=1}^s\kappa_j\nabla h_j(y)=0,\\
\lambda_i\ge 0,\;\lm_i g_i(y)=0,\;i=1,\ldots,r,\;\mbox{ and }\;\kappa_j\in \mathbb{R},\;j=1,\ldots,s
\end{array}\right\}\Longrightarrow\lambda_i=0,\;\kappa_j=0
\end{equation}
for any $y\in\O$. As we see in the subsequent sections, the imposed MMFCQ holds automatically for important classes of polynomial optimization and related problems arising in applications.}
\end{remark}\vspace*{-0.2in}

It is worth emphasizing that the majority of our applications below requires the validity of MMFCQ for the case when $\O$ in Definition~\ref{mmfcq} is chosen as the {\em argmaximum set} \eqref{Y0} at the reference point. Consider further the standard {\em Lagrangian function} for \eqref{max} with the negative sign of $\gg$ taken due to the maximization
\begin{equation}\label{La}
\mathscr{L}(x,y,\gg,\lm,\kappa):=-\gg f(x,y)+\sum_{i=1}^r\lm_ig_i(x,y)+\sum_{j=1}^s\kappa_jh_j(x,y).
\end{equation}
Observe that any $y\in M(x)$ is a minimizer of the following nonlinear program:
\begin{equation}\label{nlp}
\min-f(x,\cdot)\;\mbox{ subject to }\;g_i(x,\cdot)\le 0,\;i=1,\ldots,r,\;\mbox{ and }\;h_j(x,\cdot)=0,\,j=1,\ldots,s.
\end{equation}
Applying the classical Lagrange multiplier rule in the Fritz John form to \eqref{nlp} tells us that the set $FJ(x,y)$ of multipliers $(\gg,\lm,\kappa)\in\R_+\times\R_+^r\times\R^s$ satisfying the conditions
\begin{equation}\label{FJ0}
\gg+\sum_{i=1}^r\lm_i+\sum_{j=1}^s|\kappa_j|=1,\;\nabla_y\mathscr{L}(x,y,\gg,\lm,\kappa)=0,\;\mbox{ and }\;\lm_i g_i(x,y)=0,\;i=1,\ldots,r,
\end{equation}
is always nonempty. Given $\alpha\in[0,\infty]$, consider also the set $FJ_\alpha(x,y)$ of $(\gg,\lm,\kappa)\in\R_+\times\R_+^r\times\R^s$ satisfying
\begin{equation}\label{FJ}
\gg+\sum_{i=1}^r\lm_i+\sum_{j=1}^s|\kappa_j|\le\al,\;\nabla_y\mathscr{L}(x,y,\gg,\lm,\kappa)=0,\;\mbox{ and }\;\lm_i g_i(x,y)=0,\;i=1,\ldots,r.
\end{equation}\vspace*{-0.3in}

Using the Lagrangian description, we now show that the MMFCQ property relative to $\O=M(\ox)$ from \eqref{Y0} is {\em robust} in the general setting of Proposition~\ref{Con}(i).\vspace*{-0.1in}

\begin{proposition}{\bf (robustness of MMFCQ).}\label{MFCQS} Given $\ox\in\R^n$, assume that the sets $Y(x)$ are nonempty and uniformly bounded around $\ox$. If MMFCQ holds at $\ox$ relative to $M(\ox)$ and if $\phi$ is l.s.c.\ at $\ox$, then there exists $\delta>0$ such that MMFCQ is satisfied at any point $x\in\B_\delta(\ox)$ relative to $M(x)$.
\end{proposition}\vspace*{-0.2in}
\begin{proof}
Proposition~\ref{Lemma31} tells us that $\phi$ is continuous at $\ox$. Supposing on the contrary that MMFCQ is not robust, i.e., there are $x_k\to\ox$ as $k\to\infty$ such that MMFCQ fails at $x_k$ relative to $M(x_k)$. Take any $\xi\in\R^n$ and find $y_k\in M(x_k)$ and $(\lm_k,\kk_k)\in\R^r_+\times\R^s\setminus\{0\}$ satisfying $\nabla_y\mathscr{L}(x_k,y_k,0,\lm_k,\kappa_k)=0$, $(\lm_k)_ig_i(x_k,y_k)=0$ for $i=1,\ldots,r$, and $\nabla_x\mathscr{L}(x_k,y_k,0,\lm_k,\kappa_k)^T\xi\le 0$. Normalization gives us $\|(\lm_k,\kk_k)\|=1$ for all $k\in\mathbb{N}$. Using now the uniform boundedness of $Y(x)$ around $\ox$, we select subsequences $y_k\to\oy\in Y(\ox)$ and $(\lm_k,\kk_k)\to(\bar\lm,\bar\kk)$ with $\|(\bar\lm,\bar\kk)\|=1$. It follows from the continuity of $\phi$ at $\ox$ and from $\phi(x_k)=f(x_k,y_k)$ that $\phi(\ox)=f(\ox,\oy)$, i.e., $\oy\in M(\ox)$. Furthermore, passing to limit as $k\to\infty$ yields $\bar\lm_i g_i(\ox,\oy)=0$ for $i=1,\ldots,r$, $\nabla_y\mathscr{L}(\ox,\oy,0,\bar\lm,\bar\kappa)=0$, and $\nabla_x\mathscr{L}(\ox,\oy,0,\bar\lm,\bar\kappa)^T\xi\le 0$. Since $\xi\in\R^n$ was chosen arbitrarily, this contradicts the assumed MMFCQ at $\ox$ relative to $M(\ox)$ and thus completes the proof.\qed
\end{proof}\vspace*{-0.15in}

The next result important for its own sake plays a crucial role in deriving the extended \L ojasiewicz inequality for supremum marginal functions. It explicitly evaluates the limiting subdifferential of for such functions via the initial data and the corresponding Lagrange multiplies being significantly different from the preceding results and techniques of \cite{MN13} even for the constant mapping $Y$ in \eqref{max} considered therein. \vspace*{-0.1in}

\begin{theorem}{\bf(limiting subgradients of supremum marginal functions).}\label{lmb} Given $\ox\in\R^n$, suppose that the sets $Y(x)$ in \eqref{Y-poly} are nonempty and uniformly bounded around $\ox$, that MMFCQ holds at $\ox$ relative to $M(\ox)$, and that $\phi$ is l.s.c.\ around $\ox$. The following assertions hold:\\
{\bf (i)} There exists $\ve >0$ such that for any $x\in\mathbb{B}_\ve(\ox)$ and $v\in\partial\phi(x)$ we can find $y^{(q)}\in M(x)$ and $(\gg^{(q)},\lm^{(q)},\kk^{(q)})\in FJ_\infty(x,y^{(q)})$, $q=1,\ldots,n+1$, satisfying the conditions
\begin{equation}\label{b0}
v=\sum_{q=1}^{n+1}\nabla_x\mathscr{L}(x,y^{(q)},\gg^{(q)},\lm^{(q)},\kk^{(q)})\;\mbox{ and }\;\sum_{q=1}^{n+1}\gg^{(q)}=1,
\end{equation}
where $\mathscr{L}$ and $FJ_{\infty}$ are defined in {\rm(\ref{La})} and {\rm(\ref{FJ})}, respectively.\\
{\bf (ii)} Given $L>0$, there are positive numbers $\ve,\al$ such that for any $x\in\mathbb{B}_\ve(\ox)$ and $v\in\partial\phi(x)\cap (L\B)$ we can find $y^{(q)}\in M(x)$ and $(\gg^{(q)},\lm^{(q)},\kk^{(q)})\in FJ_\al(x,y^{(q)})$, $q=1,\ldots,n+1$, satisfying \eqref{b0}.
\end{theorem}\vspace*{-0.15in}
\begin{proof} It follows from Proposition~\ref{MFCQS} that there is $\ve>0$ such that MMFCQ holds at any point $x\in\mathbb{B}_\ve(\ox)$ relative to $M(x)$. Moreover, Proposition~\ref{Lemma31} and the l.s.c.\ assumption on $\phi$ ensure that this function is continuous on $\B_\ve(\ox)$. Let us first evaluate the proximal subdifferential \eqref{PS} and then the limiting subdifferential \eqref{ls} of the supremum marginal function \eqref{max} at each fixed $x\in\mathbb{B}_\ve(\ox)$.\vspace*{-0.1in}

{\bf Claim~1:} {\em For any proximal subgradient $v\in\partial_P\phi(x)$ we have}
\begin{equation}\label{pp}
(v,-1)\in\cone\big\{\big(\nabla_x\mathscr{L}(x,y,\gg,\lm,\kk),-\gg\big)\big|\;y\in M(x),\;(\gg,\lm,\kk)\in FJ(x,y)\big\}\subset\R^{n+1}.
\end{equation}
To verify \eqref{pp}, deduce from \eqref{PS} that there are constants $\eta>0$ and $\delta \in(0,\ve - \|x-\ox\|)$ satisfying
\begin{equation}\label{b1}
\phi(u)-\phi(x)\ge\la v,u-x\ra-\eta\|u-x\|^2\;\mbox{ whenever }\;u\in\mathbb{B}_\delta(x).
\end{equation}
The assumptions imposed on the mapping $Y$ ensure that $\phi(u)=\max\{f(u,y)|\;y\in Y(u)\}$ for all $u\in\mathbb{R}^n$, which tells us by \eqref{b1} that the pair
$(x,\phi(x))$ is a local minimizer of the following GSIP:
\begin{equation}\label{gsip1}
\min_{(u,z)}\big[z-\la v,u-x\ra+\eta\|u-x\|^2\big]\;\mbox{ subject to }\;f(u,y)-z\le 0\;\mbox{ for all }\;y\in Y(u).
\end{equation}
Applying to \eqref{gsip1} the necessary optimality condition from \cite[Theorem~1.1]{Jongen}, we find $p\in\mathbb{N}$, $y^{(q)}\in M(x)$, $(\gg^{(q)},\lm^{(q)},\kk^{(q)})\in FJ(x,y^{(q)})$, $\al\ge 0$, and $\mu^{(q)}\ge 0$ for $q=1,\ldots,p$ so that
\begin{equation}\label{b2}
\al+\sum_{q=1}^p\mu^{(q)}=1\;\mbox{ and }\;\al(-v,1)+\sum_{q=1}^p\mu^{(q)}(\nabla_x\mathscr{L}(x,y^{(q)},\gg^{(q)},\lm^{(q)},\kk^{(q)}),-\gg^{(q)})=(0,0).
\end{equation}
To justify \eqref{pp}, it remains to show that $\al>0$ in \eqref{b2}. Indeed, assuming the contrary tells us that
\begin{equation}\label{b21}
\sum_{q=1}^p\mu^{(q)}=1\;\mbox{ and }\;\sum_{q=1}^p\mu^{(q)}\big(\nabla_x\mathscr{L}(x,y^{(q)},\gg^{(q)},\lm^{(q)},\kk^{(q)}),-\gg^{(q)}\big)=(0,0)
\end{equation}
with $(\gg^{(q)},\lm^{(q)},\kk^{(q)})\in FJ(x,y^{(q)})$ for $q=1,\ldots,p$. Consider the set $I:=\big\{q\in\{1,\ldots,p\}|\;\mu^{(q)}\ne 0\big\}$, which is nonempty by $\sum_{q=1}^{p}\mu^{(q)}=1$. It follows from \eqref{b21} that $\mu^{(q)}\gg^{(q)}=0$ for $q=1,\ldots,p$. Hence $\gg^{(q)}=0$ and $(0,\lm^{(q)},\kk^{(q)})\in FJ(x,y^{(q)})$ for all $q\in I$. Combining this with \eqref{b21} yields
\[
\sum_{q\in I}\mu^{(q)}\nabla_x\mathscr{L}(x,y^{(q)},0,\lm^{(q)},\kk^{(q)})=0,
\]
which contradicts the assumed MMFCQ at $x$ relative to $M(x)$ and thus verifies \eqref{pp} for $v\in\partial_P\phi(x)$.\vspace*{-0.1in}

{\bf Claim~2:} {\em Any limiting subgradient $v\in\partial\phi(x)$ satisfies \eqref{pp}} and thus \eqref{b0}.\\[1ex]
Fix $v\in\partial\phi(x)$. Since $\phi$ is continuous on $\B_\ve(\ox)$, we find by definition \eqref{ls} sequences $x_k\to x$ and $v_k\in\partial_P\phi(x_k)$ with $v_k\to v$. Using \eqref{pp} for $v_k$ and applying the Carath\'eodory theorem to the conic convex hull in \eqref{pp} give us $y_k^{(q)}\in M(x_k)$, $\mu_k^{(q)}\in\R_+$, and $(\gg^{(q)}_k,\lm^{(q)}_k, \kk^{(q)}_k)\in FJ(x_k,y_k^{(q)})$ for $q=1,\ldots,n+1$ such that
\begin{equation}\label{b03}
v_k=\sum_{q=1}^{n+1}\mu^{(q)}_k\nabla_x\mathscr{L}(x_k,y_k^{(q)},\gg_k^{(q)},\lm_k^{(q)},\kk_k^{(q)})\;\mbox{ and }\;\sum_{q=1}^{n+1}\mu^{(q)}_k\gg_k^{(q)}=1.
\end{equation}
Let us show that the sequence of $\mu_k:=(\mu^{(1)}_k,\ldots,\mu^{(n+1)}_k)$ is bounded. Arguing by contradiction, suppose that $\sum_{q=1}^{n+1}\mu^{(q)}_k\to \infty$ and define $\hat\mu^{(q)}_k:=\mu^{(q)}_k\big[\sum_{p=1}^{n+1}\mu^{(p)}_k\big]^{-1}$ for $q=1,\ldots,n+1$. Since the convergent sequence of $v_k$ is bounded as $v_k\to v$, we deduce from \eqref{b03} that
\begin{equation}\label{b04}
\sum_{q=1}^{n+1}\hat\mu^{(q)}_k\nabla_x\mathscr{L}(x_k,y_k^{(q)},\gg_k^{(q)},\lm_k^{(q)},\kk_k^{(q)})\to 0\;\mbox{ and }\;\sum_{q=1}^{n+1}\hat\mu^{(q)}_k\gg_k^{(q)}\to 0\;\mbox{ as }\;k\to\infty.
\end{equation}
It follows from the boundedness of $y_k^{(q)}$, $\hat\mu^{(q)}_k$, and $(\gg_k^{(q)},\lm_k^{(q)},\kk_k^{(q)})$ for $q=1,\ldots,n+1$ that some subsequences of them converge to $y^{(q)}$, $\hat\mu^{(q)}$, and $(\gg^{(q)},\lm^{(q)},\kk^{(q)})$, respectively. Letting $k\to\infty$ in \eqref{b04} yields
\begin{equation}\label{b05}
\sum_{q=1}^{n+1}\hat\mu^{(q)}\nabla_x\mathscr{L}(x, y^{(q)},\gg^{(q)},\lm^{(q)},\kk^{(q)})= 0\;\mbox{ and }\;\sum_{q=1}^{n+1}\hat\mu^{(q)}\gg^{(q)}=0
\end{equation}
with $(\gg^{(q)},\lm^{(q)},\kk^{(q)})\in FJ(x,y^{(q)})$, $y^{(q)}\in M(x)$, and $\sum_{q=1}^{n+1}\hat\mu^{(q)}=1$. Similarly to Claim~1 we arrive at a contradiction with the validity of MMFCQ at $x$ relative to $M(x)$, which therefore shows that the sequence of $\mu_k=(\mu^{(1)}_k,\ldots,\mu^{(n+1)}_k)\in \R^{n+1}_+$ is bounded.\vspace*{-0.1in}

By passing to subsequences again, we get that $\mu_k\to\mu\in \R^{n+1}_+$, $y_k^{(q)}\to y^{(q)}\in M(x)$, and $(\gg_k^{(q)},\lm_k^{(q)},\kk_k^{(q)})\to(\gg^{(q)},\lm^{(q)},\kk^{(q)})\in FJ(x,y^{(q)})$ as $k\to\infty$ for $q=1,\ldots,n+1$. It follows from \eqref{b03} that
\begin{equation}\label{b06}
v=\sum_{q=1}^{n+1}\mu^{(q)}\nabla_x \mathscr{L}(x, y^{(q)},\gg^{(q)},\lm^{(q)},\kk^{(q)})\quad \mbox{and}\quad  \sum_{q=1}^{n+1}\mu^{(q)}\gg^{(q)}=1,
\end{equation}
which gives us \eqref{pp} for each $v\in\partial\phi(x)$ as stated in Claim~2.\vspace*{-0.1in}

To complete the proof of (i), define $\bar\gg^{(q)}:=\mu^{(q)}\gg^{(q)}$, $\bar\lm^{(q)}:=\mu^{(q)}\lm^{(q)}$, and $\bar\kk^{(q)}:=\mu^{(q)}\kk^{(q)}$ and then observe that the inclusion $(\gg^{(q)},\lm^{(q)},\kk^{(q)})\in FJ(x,y^{(q)})$ yields $(\bar\gg^{(q)},\bar\lm^{(q)},\bar\kk^{(q)})\in FJ_\infty(x,y^{(q)})$ and
\begin{equation}\label{bbb}
v=\sum_{q=1}^{n+1}\nabla_x \mathscr{L}(x,y^{(q)},\bar\gg^{(q)},\bar\lm^{(q)},\bar\kk^{(q)}),\quad\sum_{q=1}^{n+1}\bar\gg^{(q)}=1.
\end{equation}\vspace*{-0.3in}

To justify finally (ii), pick any $v\in\partial\phi(x)\cap (L\B)$ with $x\in\mathbb{B}_\ve(\ox)$ and find by Claim~2 and the Carath\'eodory theorem such $y^{(q)}\in M(x)$, $\mu^{(q)}\in\R_+$ and $(\gg^{(q)},\lm^{(q)},\kk^{(q)})\in FJ(x,y^{(q)})$ for $q=1,\ldots,n+1$ that \eqref{b06} holds. Taking into account that all the elements above depend on $\ve$ and proceeding as in the proof of Claim~2, we deduce that $\mu^{(q)}$, $q=1,\ldots,n+1$, are uniformly bounded when $\ve$ is sufficiently small. This means that there are numbers $\alpha,\ve>0$ such that $0\le\mu^{(q)}\le\al$ for any $\mu^{(q)}$ satisfying  \eqref{b06} whenever $x\in\B_\ve(\ox)$. Defining $\bar\gg^{(q)}:=\mu^{(q)}\gg^{(q)}$, $\bar\lm^{(q)}:=\mu^{(q)}\lm^{(q)}$, and $\bar\kk^{(q)}:=\mu^{(q)}\kk^{(q)}$ as above and then using $(\gg^{(q)},\lm^{(q)},\kk^{(q)})\in FJ(x,y^{(q)})$, we conclude that $(\bar\gg^{(q)},\bar\lm^{(q)},\bar\kk^{(q)})\in FJ_\alpha(x,y^{(q)})$ and that $v$ satisfies \eqref{bbb}. This clearly ensures the validity of \eqref{b0} and thus finishes the proof of the theorem.\qed
\end{proof}\vspace*{-0.15in}

Now we are ready to establish a subdifferential extension of the \L ojasiewicz Gradient Inequality to supremum marginal functions with polynomial data in terms of the limiting subdifferential slope \eqref{slope}.\vspace*{-0.1in}

\begin{theorem} {\bf (\L ojasiewicz subgradient inequality for polynomial systems).}\label{LojasiewiczGradientTheorem} Given $\ox\in\R^n$, suppose that all the assumptions of Theorem~{\rm\ref{lmb}} are satisfied for the marginal supremum function $\phi$ in \eqref{max}. Then there exist positive constants $c$ and $\ve$ such that we have the estimate
\begin{eqnarray}\label{MM}
{\frak m}_{\phi}(x)\ge c|\phi(x)-\phi(\bar{x})|^{1-\tau}\;\mbox{ for all }\;x\in \overline{\mathbb{B}}_{\ve}(\bar{x})\;\mbox{ with }\;\tau={\mathscr R}\big(2n+(m+r+s)(n+1),d+2\big)^{-1},
\end{eqnarray}
where the constant ${\mathscr R}$ is taken from \eqref{loj}.
\end{theorem}\vspace*{-0.15in}
\begin{proof}
It follows from Proposition~\ref{Lemma31} that $\phi$ is continuous around $\ox$. Hence there are $\nu,\ve_1>0$ such that $|\phi(x)-\phi(\ox)|\le\nu$ for all $x\in\B_{\ve_1}(\ox)$. Define $L:=\nu^{1-\tau}$ and $F\colon\mathbb{R}^n\times\mathbb{R}^m\times\R\times\mathbb{R}^r\times\mathbb{R}^s\rightarrow\mathbb{R}$
by
\begin{equation}\label{eq:oo}
F(x,y,\gg,\mu,\kappa):=-\gg f(x,y)+\sum_{i=1}^r\mu_i^2 g_i(x,y)+\sum_{j=1}^s\kappa_j h_j(x,y)
\end{equation}
and observe that $F$ is a polynomial of degree not greater than $d+2$. Define further another polynomial $P\colon\mathbb{R}^n\times\mathbb{R}^{m(n+1)}\times\R^{n}\times\mathbb{R}^{r(n+1)}\times\mathbb{R}^{s(n+1)}\rightarrow\mathbb{R}$ of degree not greater than $d+2$ by
\begin{eqnarray}\label{po}
P(x,y,\gg,\mu,\kappa):=\displaystyle\sum_{q=1}^{n}F(x,y^{(q)},\gg^{(q)},\mu^{(q)},\kappa^{(q)})+F(x,y^{(n+1)},1-\sum_{q=1}^{n}\gg^{(q)},\mu^{(n+1)},
\kappa^{(n+1)})+\phi(\ox),
\end{eqnarray}
where $x\in\R^n$, $y=(y^{(1)},\ldots,y^{(n+1)})\in\mathbb{R}^{m(n+1)}$ with $y^{(q)}\in\mathbb{R}^m$, $\gg=(\gg^{(1)},\ldots,\gg^{(n)})\in\R^n$ with $\gg^{(q)}\in\R$, $\mu=(\mu^{(1)},\ldots,\mu^{(n+1)})\in\mathbb{R}^{r(n+1)}$ with $\mu^{(q)}\in\mathbb{R}^r$, and $\kappa=(\kappa^{(1)},\ldots,\kappa^{(n+1)})\in \mathbb{R}^{s(n+1)}$ with $\kappa^{(q)}\in\mathbb{R}^s$. For each $\mu=(\mu_1,\ldots,\mu_r)\in\mathbb{R}^r$ denote $\mu^{[2]}:=(\mu_1^2,\ldots,\mu_r^2)\in\mathbb{R}^r_+$ and consider the set
\begin{eqnarray} \label{Kve}
\begin{array}{ll}
K(\ve)&\disp:=\Big\{(x,y,\gg,\mu,\kappa)\in\mathbb{R}^n\times\mathbb{R}^{m(n+1)}\times\R^{n}\times\mathbb{R}^{r(n+1)}\times \mathbb{R}^{s(n+1)}\Big|\\
&\disp \ \ \ \ \ \ x\in\overline{\mathbb{B}}_\ve(\bar{x}),\gg\in \R^{n},\sum_{q = 1}^n\gg^{(q)}\le 1,\;\gg^{(1)},\ldots,\gg^{(n)}\ge 0,\\
&\disp \ \ \ \ \ \  y=(y^{(1)},\ldots,y^{(n+1)}),\;\mu=(\mu^{(1)},\ldots,\mu^{(n+1)}),\;\kappa=(\kappa^{(1)},\ldots,\kappa^{(n+1)}),\\
&\disp \ \ \ \ \ \  y^{(q)}\in M(x),\;q=1,\ldots,n+1,\\
&\disp \ \ \ \ \ \ \left(\gg^{(q)},(\mu^{(q)})^{[2]},\kappa^{(q)}\right)\in FJ_\alpha(x,y^{(q)}),\;q=1,\ldots,n,\\
&\disp \ \ \ \ \ \ \left(1-\sum_{q=1}^{n}\gg^{(q)},(\mu^{(n+1)})^{[2]},\kappa^{(n+1)}\right)\in FJ_\alpha(x,y^{(n+1)})\Big\},\quad\ve\in(0,\ve_2),
\end{array}
\end{eqnarray}
where positive the numbers $\alpha$ and $\ve_2\le\ve_1$ depend on $L$ and are taken from Theorem~\ref{lmb}(ii). It can be verified that the set-valued mapping $K(\cdot)$ is o.s.c.\ at $\ve=0$. Furthermore, for any $(x,y,\gg,\mu,\kappa)\in K(0)$ we get from \eqref{FJ} and \eqref{Kve} the equalities
\[\begin{array}{ll}
P(x,y,\gg,\mu,\kappa)&\disp=-\sum_{q=1}^{n}\gg^{(q)}f(x,y^{(q)})-\Big((1-\sum_{q=1}^{n}\gg^{(q)})f(x,y^{(n+1)}\Big)+\phi(\ox)\\
&\disp=-\sum_{q=1}^{n}\gg^{(q)}\phi(\ox)-\Big(1-\sum_{q=1}^{n}\gg^{(q)}\Big)\phi(\ox)+\phi(\ox)=-\phi(\ox)+\phi(\ox)=0,
\end{array}
\]
where the second one follows from $y^{(q)}\in M(x)=M(\ox)$ as $q=1,\ldots,n+1$. Applying the \L ojasiewicz Gradient Inequality \eqref{loj} to the polynomial $P$ in \eqref{po} gives us positive constants $\eta$ and $c_\eta$ such that
\begin{eqnarray}\label{b5}
\|\nabla P(x',y',\gg',\mu',\kappa')\|\ge c_\eta\,|P(x',y',\gg',\mu',\kappa')|^{1-\tau}\;\mbox{ if }\;\|(x',y',\gg',\mu',\kappa')-(x,y,\gg,\mu,\kappa)\|<\eta
\end{eqnarray}
with $\tau:={\mathscr R}(2n+(m+r+s)(n+1),d+2)$. It is easy to see that $K(0)$ is a compact set, and hence we can find $p\in\mathbb{N}$, $(x^{(\nu)},y^{(\nu)},\gamma^{(\nu)},\mu^{(\nu)},\kappa^{(\nu)})\in K(0)$, and $\eta_{\nu}>0$ for $\nu=1,\ldots,p$ so that
\[
K(0)\subset\bigcup_{\nu=1}^p\mathbb{B}_{\eta_\nu}(x^{(\nu)},y^{(\nu)},\gamma^{(\nu)},\mu^{(\nu)},\kappa^{(\nu)}).
\]
This allows us to find a positive number $\delta$ ensuring the inclusion
\[
K(0)+\delta\mathbb{B}\subset\bigcup_{\nu=1}^p\mathbb{B}_{\eta_\nu}\big(x^{(\nu)},y^{(\nu)},\gamma^{(\nu)},\mu^{(\nu)},\kappa^{(\nu)}\big).
\]
Employing the o.s.c.\ property of $K(\cdot)$ at $\ve=0$, we find $\ve_3\in(0,\ve_2)$ so that
\[
K(\ve)\subset K(0)+\delta\mathbb{B}\subset\bigcup_{k=1}^p\mathbb{B}_{\eta_k}\big(x^{(\nu)},y^{(\nu)},\gamma^{(\nu)},\mu^{(\nu)},\kappa^{(\nu)}\big)
\]
for all $\ve\in[0,\ve_3)$. Denoting $c:=\min\{1,c_{\eta_1},\ldots,c_{\eta_p}\}$, deduce from \eqref{b5} that
\begin{eqnarray}\label{eq:didio}
\|\nabla P(x,y,\gg,\mu,\kappa)\|\ge c\,|P(x,y,\gg,\mu,\kappa)|^{1-\tau}\;\mbox{ if }\;(x,y,\gg,\mu,\kappa)\in K(\ve).
\end{eqnarray}
Pick now $\ve\in(0,\ve_3)$ and $x\in\overline{\mathbb{B}}_\ve(\bar{x})$. If $m_\phi(x)\ge L$, then \eqref{MM} holds by the choices of $L$ and $c$. If $m_\phi(x)<L$, then we take any $v\in\partial\phi(x)\cap (L\B)\ne\emp$ and deduce from Theorem~\ref{lmb}(ii) that
there are $y^{(q)}\in M(x)$ and $(\gg^{(q)},\lm^{(q)},\kk^{(q)})\in FJ_\alpha(x,y^{(q)})$, $q=1,\ldots,n+1$, satisfying
\begin{equation}\label{b6}
v=\sum_{q=1}^{n+1}\nabla_x\mathscr{L}(x,y^{(q)},\gg^{(q)},\lm^{(q)},\kk^{(q)})\;\mbox{ and }\;\sum_{q=1}^{n+1}\gg^{(q)}=1.
\end{equation}
Define next $\mu^{(q)}_i:=\sqrt{\lambda^{(q)}_i}$, $i=1,\ldots,r$, and observe the inclusions
\[
(\gg^{(q)},[\mu^{(q)}]^2,\kk^{(q)})\in FJ_{\alpha}(x,y^{(q)}),\quad q=1,\ldots,n+1.
\]
It follows from \eqref{b6} and the definition of $K(\ve)$ in \eqref{Kve} that $(x,y,\gg,\mu,\kappa)\in K({\ve})$, and thus (\ref{eq:didio}) holds for this quintuple. Differentiating \eqref{po} with taking \eqref{b6} into account implies the equalities
\begin{eqnarray*}\begin{array}{ll}
\nabla_x P(x, {y},\gg,\mu,\kappa)&\disp=\sum_{q=1}^{n+1}\nabla_x\mathscr{L}(x,y^{(q)},\gg^{(q)},\lm^{(q)},\kk^{(q)})=v,\\
\nabla_y P(x,y,\gg,\mu,\kappa)&\disp=\Big(\nabla_{y^{(1)}}\mathscr{L}(x,y^{(1)},\gg^{(1)},\lm^{(1)},\kk^{(1)}),\ldots,\nabla_{y^{(n+1)}}
\mathscr{L}(x,y^{(n+1)},\gg^{(n+1)},\lm^{(n+1)},\kk^{(n+1)})\Big)=0,\\
\nabla_{\gg} P(x,y,\gg,\mu,\kappa)&\disp=\big(f(x,y^{(1)})-f(x,y^{(n+1)}),\ldots,f(x,y^{(n)})-f(x,y^{(n+1)})\big)=0,\\
\nabla_\mu P(x,y,\gg,\mu,\kappa)&\disp=\Big(2\mu^{(q)}_ig_i(x,y^{(q)})\Big)_{i,q}=0,\quad\mbox{and}\quad \nabla_{\kappa}P(x,y,\gg,\mu,\kappa)=\Big(h_j(x,y^{(q)})\Big)_{j,q}=0,
\end{array}
\end{eqnarray*}
which ensure by \eqref{eq:didio}, \eqref{b6}, and \eqref{po} the relationships
\begin{eqnarray*}
\|v\|=\|\nabla P(x,y,\gamma,\mu,\kappa)\|\ge c\,|P(x, {y},\gamma,\mu,\kappa)|^{1-\tau}=c\,|\phi(x)-\phi(\bar{x})|^{1-\tau}\;\mbox{ for any  }\;v\in\partial\phi(x)\cap(L\B).
\end{eqnarray*}
Since $m_\phi(x)<L$, this clearly verifies \eqref{MM} and completes the proof of the theorem.\qed
\end{proof}\vspace*{-0.15in}
\begin{remark}{\bf(sharper estimate for equality polynomial systems).}\label{remark:1} A close look at the proof of Theorem~\ref{LojasiewiczGradientTheorem} reveals that the exponent $\tau$ in \eqref{MM} can sharpen to $\tau={\mathscr R}(2n+(m+s)(n+1),d+1)^{-1}$ in the case when the polynomial system in \eqref{Y-poly} is given by only the equalities $Y(x)=\{y\in\mathbb{R}^m|\;h_j(x,y)=0,\;j=1,\ldots,s\}$. This is due to the fact that the polynomials $F$ in (\ref{eq:oo}) and $P$ in \eqref{po} are of degrees at most $d+1$ instead of $d+2$ as in the general case.
\end{remark}
\vspace*{-0.3in}

\section{Error Bounds for Parametric Polynomial Systems}\vspace*{-0.1in}
\setcounter{equation}{0}

In this section we derive H\"{o}lder error bounds with explicit exponents for general polynomial systems of type \eqref{eq:S}, where the underlying state-dependent sets of parameters $Y(x)$ admit the polynomial description \eqref{Y-poly}. We also present important specifications of these estimates in some particular cases. The results obtained are far-going extensions of those established in \cite{LMP2014} in the case of finite and constant sets $Y(x)\equiv\O$.
To proceed, we need the following lemma on error bounds for locally Lipschitz functions taken from \cite[Corollary~2]{Thera}.\vspace*{-0.1in}

\begin{lemma}{\bf (sufficient condition for error bounds of Lipschitzian functions).}\label{lemma:thera_errorbound} Let $f\colon\mathbb{R}^n\rightarrow\mathbb{R}$ be continuous around $\overline{x}\in\bd S_f$, where $S_f:=\{x|\;f(x)\le 0\}$. Assume that there are numbers $c,\ve>0$ such that ${\frak m}_f(x)\ge c\,|f(x)|^{1-\tau}$ for all $x\in\R^n$ with $\|x-\overline{x}\|\le\ve$ and $x\notin S_f$. Then we have
\[
{\rm dist}(x,S_f)\le\frac{1}{\tau c}\big[f(x)\big]_+^{\tau}\;\mbox{ whenever }\;\|x-\overline{x}\|\le\frac{\ve}{2}.
\]
\end{lemma}\vspace*{-0.15in}

We are now ready to establish the main result of this section. To avoid some technical complexity on dealing with the l.s.c.\ property of supremum marginal functions as discussed in Section~3, we assume in what follows that the mapping $Y$ in \eqref{Y-poly} is {\em inner semicontinuous} around the point in question, which is always the case when $Y(x)\equiv \O$ for some fixed set $\O\subset\R^m$; see, e.g., \eqref{Omega0} below.\vspace*{-0.1in}

\begin{theorem}{\bf (H\"{o}lder error bounds for parametric polynomial systems).}\label{th:local}  Let $f_l,g_i,h_j\colon\mathbb{R}^n\times\mathbb{R}^m\rightarrow\mathbb{R}$ with $l=1,\ldots,L$, $i=1,\ldots,r$, and $j=1,\ldots,s$ be polynomial functions of degrees at most $d$. Suppose that the mapping $Y$ from \eqref{Y-poly} is i.s.c.\ with nonempty and uniformly bounded values $Y(x)$ around a given point $\ox\in\R^n$ and that MMFCQ from Definition~{\rm\ref{mmfcq}} holds at this point relative to the argmaximum set
\begin{equation}\label{argmax}
M(\ox):=\big\{y\in Y(\bar x)\big|\;\displaystyle\max_{y\in Y(x)\atop 1\le l\le L}f_l(\ox,y)=f_l(\ox,y)\;\mbox{ for some }\;l=1,\ldots,L\big\}.
\end{equation}
Then there are numbers $c,\ve>0$ such that we have the local error bound
\begin{equation}\label{ib}
{\rm dist}(x,S)\le c\,\bigg(\Big[\max_{y\in Y(x)\atop 1\le l\le L}f_l(x,y)\Big]_{+}\bigg)^{\tau}\;\mbox{ for all }\;x\in\R^n\;\mbox{ with }\;\|x-\bar x\|\le \ve,
\end{equation}
where $S$ is taken from \eqref{eq:S}, and where the H\"older exponent $\tau$ in \eqref{ib} is explicitly calculated by
\begin{equation}\label{tau1}
\tau=\left\{\begin{array}{ll}
{\mathscr R}(2n+(m+r+s)(n+1),d+2)^{-1},&L=1,\\
{\mathscr R}(2n+(m +r+s+2)(n+1),d+L+1)^{-1},&L\ge 2,
\end{array}\right.
\end{equation}
with the constant ${\mathscr R}$ defined in \eqref{loj}. If we suppose in addition that MMFCQ holds at every $x\in\R^n$ relative to the whole set of parameters $Y(x)$ and that $Y$ are i.s.c on $\R^n$, then for any compact set $K\subset\mathbb{R}^n$ there is a number $\bar c>0$ such that we have the following global error bound with the same exponent $\tau$:
\begin{equation}\label{ib4}
{\rm dist}(x,S)\le\bar c\, \bigg(\Big[\max_{y\in Y(x)\atop 1\le l\le L}f_l(x,y)\Big]_{+}\bigg)^{\tau}\;\mbox{ whenever }\;x\in K.
\end{equation}
\end{theorem}\vspace*{-0.15in}
\begin{proof}
Consider first the case of $L=1$, where $\phi(x)=\sup_{y\in Y(x)}f_1(x,y)$ is continuous around $\ox$ by Propositions~\ref{Lemma31} and \ref{Con}(i). Theorem~\ref{LojasiewiczGradientTheorem} gives us numbers $c_0,\ve_0$ such that for all $x\in\overline{\mathbb{B}}_{\ve_0}(\bar{x})$ we have the estimate
$$
{\frak m}_{\phi}(x)\ge c_0\,|\phi(x)-\phi(\bar{x})|^{1-\tau}\;\mbox{ with }\;\tau={\mathscr R}\big(2n+(m+r+s)(n+1,d+2\big)^{-1}.
$$
Let us justify the existence of $c,\ve>0$ for which \eqref{ib} holds with $\tau$ from \eqref{tau1} as $L=1$. If $\bar x\in{\rm bd}\,S$, i.e., $\phi(\ox)=0$, we get by Lemma~\ref{lemma:thera_errorbound} that \eqref{tau1} is fulfilled with $c=\frac{1}{\tau c_0}$ and $\ve=\ve_0/2$. Since this inequality is automatic for $\bar{x}\in{\rm int}\,S$, it remains to examine the case of $\bar x\notin S$. It follows from $\bar x\not\in S=\{x\;|\;\phi(x)\le 0\}$ that $\phi(\bar x)>0$ and so $(\phi(\bar x))^\tau>0$. By the continuity of $\phi$ and $d(\cdot,S)$ we find $M_1,M_2,\ve>0$ such that
\[
(\phi(x))^\tau\ge M_1\;\mbox{ and }\;{\rm dist}(x,S)\le M_2\;\mbox{ for all }\;x\in\overline{\mathbb{B}}_{\ve}(\bar{x}),
\]
and hence ${\rm dist}(x,S)\le M_2\le\frac{M_2}{M^\tau_1}\phi(x)^\tau$ for all $x$ from this ball. This verifies \eqref{ib} for $L=1$.\vspace*{-0.1in}

Next we consider the case of $L\ge 2$ and define the {\em Lagrange interpolation polynomials} $\gamma_{l}\colon\mathbb{R}\rightarrow\mathbb{R}$ by
\[
\gamma_l(t):=\prod_{j\ne l,j\in\{1,\ldots,L\}}\frac{t-j}{l-j},\quad l=1,\ldots,L,
\]
for which $\gamma_l(l)=1$ as $l=1,\ldots,L$ and $\gamma_l(j)=0$ as $j\in\{1,\ldots,L\}$ with $j\ne l$. Define the functions
\begin{eqnarray*}
\begin{array}{ll}
\Tilde{f}(x,y,t):=\disp\sum_{l=1}^{L}\gamma_l(t)f_l(x,y),\quad\Tilde g_i(x,y,t):=g_i(x,y),\;i=1,\ldots,r,\\
\Tilde h_j(x,y,t):=h_j(x,y),\;j=1,\ldots,s,\;\mbox {and }\;\Tilde h_{s+1}(x,y,t):=\disp\prod_{l=1}^{L}(t-l)
\end{array}
\end{eqnarray*}
and observe that $\Tilde{f}$ is a polynomial on $\mathbb{R}^n\times\mathbb{R}^{m+1}$ with degree not greater than $d+L-1$. Form now the modified constraint and argmaximum sets by, respectively,
\begin{eqnarray*}
\Tilde Y(x):=\big\{(y,t)\in\R^m\times\R\big|\;\Tilde g_i(x,y,t)\le 0,\;i=1,\ldots,r;\;\Tilde h_j(x,y,t)=0,\;j=1,\ldots,s+1\big\},
\end{eqnarray*}
\begin{eqnarray*}
\Tilde M(x):=\big\{(y,t)\in\Tilde Y(x)\big|\;\Tilde f(x,y,t)=\max_{(y,t)\in\Tilde Y(x)}\Tilde{f}(x,y,t)\big\}
\end{eqnarray*}
while noting that $\Tilde Y(x)=Y(x)\times\{1,\ldots,L\}$ for $Y(x)$ from \eqref{Y-poly} and that
\[
\max_{y\in Y(x)\atop 1\le l\le L}f_l(x,y)=\max_{(y,t)\in\Tilde Y(x)}\Tilde{f}(x,y,t),\quad S=\big\{x\in\R^n\big|\;\Tilde f(x,y,t)\le 0\;\mbox{ for all } \;(y,t)\in\Tilde Y(x)\big\}.
\]
It can be directly checked that the mapping $\Tilde Y$ is i.s.c.\ around $\ox$ with the nonempty and uniformly bounded values and that the MMFCQ property for the modified system holds at $\ox$ relative to $\Tilde M(\ox)$. Applying the result for the case of $L=1$  obtained above to the system $\Tilde{f},\Tilde g_i,\Tilde h_j$ with $(n,m,r,s,d)$ replaced by $(n,m+1,r,s+1,d+L-1)$, we arrive at \eqref{ib} with $\tau$ calculated in \eqref{tau1} for $L\ge 2$. The last statement of the theorem \eqref{ib4} follows from \eqref{ib} by employing the standard compactness arguments.\qed
\end{proof}\vspace*{-0.15in}

\begin{remark}{\bf (sharper error bounds for equality polynomial systems).}\label{remark:21} Based on the observation in Remark~\ref{remark:1} and following the proof of Theorem~\ref{th:local}, we see that the the exponent estimate $\tau$ in \eqref{ib} can be sharpen to $\tau={\mathscr R}(2n+(m+r+s+2\min\{L-1,1\})(n+1),d+L)$ if the polynomial system $Y(x)$ in \eqref{Y-poly} contains only the equality constraints.
\end{remark}\vspace*{-0.15in}

In the rest of this section we derive specifications of Theorem~\ref{th:local} in some particular classes of parametric polynomial systems and discuss the possibility of sharpening H\"older exponents in certain situations.\vspace*{-0.2in}

\subsection*{Error bounds for infinite polynomial inequality systems}\vspace*{-0.1in}

Let us examine the case of \eqref{eq:S} where the sets $Y(x)$ are {\em independent} of $x$ being reduced to the constant set
\begin{equation}\label{Omega0}
\Omega:=\big\{y\in\R^m\big|\;g_i(y)\le 0,\;i=1,\ldots,r;\;h_j(y)=0,\;j=1,\ldots,s\big\}
\end{equation}
described by polynomials with degrees at most $d$. In particular, such a case covers important models known as {\em polynomial matrix inequalities}; see below. Observe also that $\O$ from \eqref{Omega0} can be treated as an infinite index set in problems of {\em semi-infinite programming} (SIP); see, e.g., \cite{GL01}. However, semialgebraic structures of index sets have not been specifically investigated in the usual framework of SIP.\vspace*{-0.1in}

Following the traditional terminology, we say that $\O$ is {\em regular} if MFCQ for \eqref{Omega0} holds on $\O$. Let us present examples of two important classes of regular sets described by polynomials.\vspace*{-0.1in}

$\bullet$ {\bf Regularity of discrete sets}. Considering the discrete set $\Omega:=\{1,\ldots,p\}$, $p\in\mathbb{N}$, we can write it as
$$
\Omega=\big\{y\in \mathbb{R}\big|\;h(y)=0\big\}\;\mbox{ with }\;h(y):=(y-1)(y-2)\cdots(y-p).
$$
Since we obviously have $\nabla h(y)\ne 0$ for all $y\in\Omega$, this set is regular.\vspace*{-0.1in}

$\bullet$ {\bf Regularity of algebraic manifolds}. It follows directly from the definition that the algebraic manifold
$$
\O:=\big\{y\in\R^m\big|\;h_1(y)=0,\ldots,h_s(y)=0\big\},\quad s\le m,
$$
is regular provided that $\mathrm{rank}\big(\nabla h_1(y),\ldots,\nabla h_s(y)\big)=s$ for all $y\in\O$.\vspace*{-0.05in}

As an immediate consequence of Theorem~\ref{th:local}, we arrive at the following extension of the results in \cite{LMP2014} from finite to to infinite polynomial inequality systems.\vspace*{-0.1in}

\begin{corollary}{\bf (H\"older error bounds for infinite polynomial systems).}\label{cor:local} Let $f_l\colon\mathbb{R}^n\times\mathbb{R}^m\to\mathbb{R}$, $l=1,\ldots,L$, be polynomials of degrees at most $d$, let $\O\ne\emp$  given  in \eqref{Omega0} be bounded and regular, and let $S$ be defined in \eqref{eq:S} with $Y(x)\equiv\O$. Then for any $\ox\in\R^n$ there exist constants $c,\ve>0$ such that
\begin{equation}\label{ib99}
{\rm dist}(x,S)\le c\,\bigg(\Big[\max_{y \in\Omega\atop 1\le l\le L}f_l(x,y)\Big]_{+}\bigg)^{\tau}\;\mbox{ whenever }\;\|x-\bar x\|\le\ve,
\end{equation}
where $\tau={\mathscr R}(2n+(m+r+s+2\min\{L-1,1\})(n+1),d+L+1)^{-1}$ with the constant ${\mathscr R}$ taken from \eqref{loj}. Furthermore, the error bound \eqref{ib99} holds globally with some $c>0$ on any compact set $K\subset\R^n$.
\end{corollary}\vspace*{-0.2in}

\begin{remark}{\bf (semialgebracity of index sets is essential for H\"older error bounds).}\label{remark:6} As follows from Example~\ref{ex:1}, the error bound (\ref{ib99}) can fail if $\O$ is not semialgebraic. Indeed, this example shows that the H\"older error bound does not hold whenever $\tau\in(0,1)$ for $f$ and $\O$ therein. It remains to check that the set $\O$ from \eqref{exO} is not semialgebraic. Supposing the contrary implies by the Tarski-Seidenberg Theorem that the function $x\mapsto\max_{y\in\Omega}f(x,y)$ is semialgebraic on $[-0.5,0.5]$, which is not the case.
\end{remark}\vspace*{-0.15in}

Next we present a consequence of Corollary~\ref{cor:local} to polynomial matrix inequalities. Recall that an $(m\times m)$ {\em polynomial matrix inequality} of  $n$ variables with degree $d$ is represented in the form $P(x)\preceq 0$, where the notation $M\preceq 0$ means that $-M$ is a positive semidefinite matrix, and where $P\colon\mathbb{R}^n\rightarrow S^m$ is a mapping such that each $(i,j)$th element of its image $\big(P(x)\big)_{ij}$, $1\le i,j\le m$, is a real polynomial of degree at most $d$. We refer the reader to, e.g., \cite{Las_Henrion,Shuzhong} for various applications of such systems.\vspace*{-0.1in}

\begin{corollary}{\bf (error bounds for polynomial matrix inequalities without regularity assumptions).}\label{cor:PMI} Fix $d\in\mathbb{N}$ and consider the solution set
$$
S_{PMI}:=\big\{x\big|\;P(x)\preceq 0\big\}
$$
described by the above polynomial matrix inequality. Then for a compact set $K\subset\R^n$ there is $c>0$ such that
\[
{\rm dist}(x,S_{PMI})\le c\,\left(\Big[\lambda_{\rm max}\big(P(x)\big)\Big]_+\right)^{\tau}\;\mbox{ whenever }\;x \in K,
\]
where $\tau={\mathscr R}(2n+(m+1)(n+1),d+3)^{-1}$ and $\lambda_{\max}$ denotes the maximum eigenvalue of a symmetric matrix.
\end{corollary}\vspace*{-0.1in}
\begin{proof}
Define $f(x,y):=\sum_{i,j=1}^m y_i\left(P(x)\right)y_j$, which is a real polynomial of degree $d+2$ on $\mathbb{R}^{n+m}$. We have
\begin{equation}\label{PMI1}
P(x)\preceq 0\Longleftrightarrow\Big[f(x,y)\le 0\;\mbox{ for all }\;y\in\O:=\big\{y\in\mathbb{R}^m\big|\;\|y\|^2=1\big\}\Big].
\end{equation}
Note that the set $\O\ne\emp$ in \eqref{PMI1} is bounded and regular. Applying now Corollary~\ref{cor:local} to the inequality system on the right-hand side of \eqref{PMI1} and replacing $d$ by $d+2$ with taking into account Remark~\ref{remark:1} ensure that for any compact set $K$ we arrive at the relationships
\[
{\rm dist}(x,S_{PMI})\le c\,\left(\Big[\sup_{\|y\|^2=1}f(x,y)\Big]_+\right)^{\tau}=c\,\left(\Big[\lambda_{\rm max}\big(P(x)\big)\Big]_+\right)^{\tau}
\]
with some constant $c>0$. This completes the proof of the corollary.\qed
\end{proof}\vspace*{-0.2in}

\subsection*{Discussion on sharpness of exponent estimations}\vspace*{-0.1in}

Finally in this section, we discuss the possibility to sharpen H\"older exponents in the established error bounds for general infinite polynomial systems in some particular settings. Let us first present an example of a polynomial matrix inequality, where the exponent obtained in Corollary~\ref{cor:PMI} is not the sharpest one.\vspace*{-0.1in}

\begin{example} {\bf (sharper error bound for matrix inequalities).}\label{tighter} Let $d$ be an even number. For any $x=(x_1,\ldots,x_n)$ define $A_1(x):=x_1^d$ and then, for any $i=2,\ldots,n$,
$$
A_i(x):=\left(\begin{array}{ccc}
-1&x_{i}^d\\x_{i}^d&-x_{i-1}
\end{array}\right)=\left(\begin{array}{cc}
-1&0\\
0&0
\end{array}\right)+x_{i-1}\left(\begin{array}{cc}
0&0\\
0&-1
\end{array}\right)+x_{i}^d \left(\begin{array}{cc}
0&1\\
1&0
\end{array}\right).
$$
Consider the following $(2n-1)\times (2n-1)$ polynomial matrix inequality system:
\[P(x):=\left(\begin{array}{cccc}
A_1(x)&0&0&0\\
0&A_2(x)&0&0\\
\vdots&0&\ddots&0\\
0&0&\ldots&A_{n}(x)
\end{array}\right)\preceq 0
\]
and get by direct calculations that $S_{PMI}:=\{x\in\mathbb{R}^n|\;P(x)\preceq 0\}=\{0\}$ and
\begin{eqnarray*}
\lambda_{\rm max}\big(P(x)\big)&=&\max\big\{\lambda_{\rm max}\big(A_i(x)\big)\big|\;1\le i\le n-1\big\}\\
&=&\max\Big\{\max\Big\{\frac{-(1+x_{i-1})}{2}+\sqrt{\frac{(1-x_{i-1})^2}{4}+x_{i}^{2d}},\;x_{1}^d\Big\}\Big|\;2\le i\le n\Big\} \\
&=&\max\left\{\max\left\{\frac{x_i^{2d}-x_{i-1}}{\frac{(1+x_{i-1})}{2}+\sqrt{\frac{(1-x_{i-1})^2}{4}+x_{i}^{2d}}},\;x_{1}^d\right\}\Big|\;2\le i\le n\right\}.
\end{eqnarray*}
Observe further that for the parametric family $x(t):=(t^{(2d)^{n-1}},\ldots,t^{2d},t)$, $t\in(0,1)$, we have $d(x(t),S_{PMI})=O(t)$ and
$\lambda_{\max}\big(P(x(t))\big)=t^{d(2d)^{n-1}}$. Thus the exponent of the error bound for $\{x\in\mathbb{R}^n|\;P(x)\preceq 0\}$ at the origin is at most $\tau=\frac{1}{d(2d)^{n-1}}$. Moreover, there is a constant $c_0>0$ such that for all $x$ around $\ox=0$ we have the relationships
\begin{eqnarray*}
\lambda_{\rm max}\big(P(x)\big)
&=&\max\left\{\max\left\{\frac{x_i^{2d}-x_{i-1}}{\frac{(1+x_{i-1})}{2}+\sqrt{\frac{(1-x_{i-1})^2}{4}+x_{i}^{2d}}},\;x_{1}^d\right\}\Big|\;2\le i\le n\right\}\\
&=&\max\left\{\max\left\{\left[\frac{x_i^{2d}-x_{i-1}}{\frac{(1+x_{i-1})}{2}+\sqrt{\frac{(1-x_{i-1})^2}{4}+x_{i}^{2d}}}\right]_+,\;x_{1}^d\right\}\Big|\;2\le i\le n\right\}\\
 &\ge& \max\Big\{\max\Big\{c_0\,\big[x_i^{2d}-x_{i-1}\big]_+,\;x_{1}^d\Big\}\Big|\;2\le i\le n\Big\},
\end{eqnarray*}
where the second equality is due to $x_1^d\ge 0$ while the last inequality follows as the denominator in the third equality approaches one as $x\rightarrow 0$. Considering next the finite convex polynomial system
\[
x_{1}^d\le 0\;\mbox{ and }\;c_0\,(x_i^{2d}-x_{i-1})\le 0\;\mbox{ for }\;i=2,\ldots,n,
\]
note that its solution set is the same as $S_{PMI}=\{0\}$. It is known \cite{Borwein} that for any finite convex polynomial system, where all the involved polynomials have degree at most $d_0$ on $\mathbb{R}^n$, a local H\"{o}lder error bound holds with the exponent $\frac{2}{(2d_0-1)^n+1}$. Hence there are numbers $c,\delta>0$ such that
$$
d(x,S_{PMI})\le c\max\big\{\max\big\{c_0\,[x_i^{2d}-x_{i-1}]_+,\;x_{1}^d\big\}\big|\;2\le i\le n\big\}^{\frac{2}{(4d-1)^n+1}}\le c\,\lambda_{\rm max}\big(P(x)\big)^{\frac{2}{(4d-1)^n+1}},\quad\|x\|\le\dd.
$$
Thus the {\em true exponent} of the H\"older error bound at the origin in this example lies within the interval $[{\frac{2}{(4d-1)^n+1}},\frac{1}{d (2d)^{n-1}}]$.  On the other hand, the exponent value calculated in Corollary~\ref{cor:PMI} is
$$
\tau=\disp\frac{1}{{\mathscr R}(2n^2+4n,d+3)}=\frac{1}{(d+3)(3d+6)^{2n^2+4n-1}},
$$
which therefore gives us only a lower estimate of the true H\"older exponent in the local error bound.
\end{example}\vspace*{-0.15in}

We now list two specific structures of polynomial systems, where error bounds can be {\em sharpened}:\vspace*{-0.1in}

$\bullet$ Let $(x,y)\mapsto f_l(x,y)$ be a {\em concave polynomial} for all $l=1,\ldots, L$, $\Omega$ be a compact convex set, and there exist $x_0\in\mathbb{R}^n$ with $\displaystyle\max_{y\in\Omega,\atop 1\le l\le L}f_l(x_0,y)<0$. It can be directly verified in this case that $\displaystyle x \mapsto\max_{y\in\Omega,\atop 1\le l\le L}f_l(x,y)$ is a continuous concave function. Note that a {\em Lipschitz local error bound }(i.e., with the best possible exponent $\tau=1$) holds for a concave function satisfying the imposed strict feasibility condition; see \cite[Theorem~3]{Ngai}.\vspace*{-0.1in}

$\bullet$ Let $x\mapsto f_l(x,y)$ be a {\em convex polynomial} for each $y\in\mathbb{R}^m$ and $l=1,\ldots,L$, and let $\Omega$ be a convex polytope in $\mathbb{R}^m$. Since in this case the number of extreme points of $\Omega$ is finite and is at most $m+1$, we get the representation $\Omega={\rm co}\{y^{(1)},\ldots,y^{(m+1)}\}$. Noting that a convex function attains its maximum over a convex polytope at its extreme points, it follows that
$
\disp\max_{y\in\Omega,\atop 1\le l\le L}f(x,y)=\max_{1\le i\le m+1}f_l(x,y^{(i)}).
$
Thus the corresponding explicit error bounds can be simplified and sharpened in this case by using the results from our preceding paper \cite{LMP2014} for finite index sets $\O$.\vspace*{-0.2in}

\section{Higher-Order Stability Analysis for Optimization Problems}\vspace*{-0.1in}
\setcounter{equation}{0}

In this section we develop some applications of the error bound results from Section~4 to higher-order stability analysis for three important classes in {\em parametric polynomial optimization} dealing with {\em infinite} sets of parameters $Y(x)$ from \eqref{Y-poly}. In contrast to the majority of publications in parametric optimization, our main emphasis is on {\em higher-order stability} of set-valued solution maps with {\em explicit calculating} the corresponding exponents. Let us start with problems of generalized semi-infinite programming.\vspace*{-0.2in}

\subsection{Generalized Semi-Infinite Programs with Polynomial Data}\vspace*{-0.1in}

Consider the following {\em polynomial generalized semi-infinite program}:
\begin{eqnarray*}
(GSIP) &\qquad\mbox{minimize}_{x\in\R^n} &p_0(x),\\
&\qquad\mbox{ subject to }&f_l(x,y)\le 0\;\mbox{ for all }\;y\in Y(x),\;l=1,\ldots,L,
\end{eqnarray*}
where the sets $Y(x)\subset\mathbb{R}^m$ are taken from \eqref{Y-poly}, and where all the functions $p_0$, $f_l$, $g_i$, and $h_j$ are real polynomials with degrees at most $d$. Along with problem $(GSIP)$ we consider also its {\em perturbation}
\begin{eqnarray*}
(GSIP_u) &\qquad\mbox{minimize}_{x\in\R^n}& p_0(x,u),\\
&\qquad\mbox{ subject to }&f_l(x,y)\le 0\;\mbox{ for all }\;y\in Y(x),\;l=1,\ldots,L,
\end{eqnarray*}
where the cost function $p(x,u)$ with $p(x,0)=p_0(x)$ depends on the perturbation parameter $u\in\R^q$ being polynomial of degree at most $d$ in $x$ and locally Lipschitzian in $u$. Note that the case of $p(x,u)=p_0(x)-\la u,x\ra$ and $Y(x)\equiv\{0\}$ in $(GSIP_u)$ corresponds to the so-called {\em tilt perturbation} introduced by Poliquin and Rockafellar in \cite{pr}, where the focus is on {\em single-valuedness} and {\em Lipschitz continuity} of the optimal solution map {\em around} the point in question, which is known as {\em tilt stability}. In \cite{gm,mr} the reader can find recent second-order characterizations of tilt stability in the aforementioned setting of nonlinear programming with ${\cal C}^2$-smooth data and related issues. The {\em quantitative higher-order} stability for {\em polynomial} GSIPs established in the next theorem is significantly different: we justify a certain {\em H\"older continuity} of the {\em set-valued} solution map {\em at} the reference point with calculating the corresponding H\"older exponent.\vspace*{-0.1in}

\begin{theorem}{\bf (quantitative higher-order stability for polynomial GSIPs).}\label{th:stability_SIP} Let the feasible set
\begin{equation}\label{G}
G:=\big\{x\in\R^n\big|\;f_l(x,y)\le 0\;\mbox{ for all }\;y\in Y(x),\;l=1,\ldots,L\big\}
\end{equation}
be nonempty and compact. Assume that the mapping $Y$ in \eqref{G} is i.s.c.\ on $\R^n$ with nonempty and uniformly bounded values and that MMFCQ holds on $\R^n$ relative to $Y(x)$. Denote by $S(u)$ the set of optimal solutions to the perturbed problem $(GSIP_u)$. Then for any fixed $\bar u\in\mathbb{R}^q$ there are constants $c,\ve>0$ such that
\begin{equation}\label{Su}
S(u)\subset S(\bar u)+c\,\|u-\bar u\|^\tau\;\mbox{ if }\;\|u-\ou\|\le\ve\;\mbox{ with }\;\tau={\mathscr R}(2n+(m+r+s+2)(n+1),d+L+2)^{-1}.
\end{equation}
\end{theorem}\vspace*{-0.1in}
\begin{proof} Consider the {\em optimal value function} in $(GSIP_u)$ given by
$$
p_*(u):=\inf\big\{p(x,u)\big|\;f_l(x,y)\le 0\;\mbox{ whenever }\;y\in Y(x),\;l=1,\ldots,L\big\},\quad u\in\R^q,
$$
and observe that the solution map $S(u)$ to $(GSIP_u)$ is represented as
\begin{equation}\label{phi-sip0}
S(u)=\big\{x\in\mathbb{R}^n\big|\;f_l(x,y)\le 0\;\mbox{ for all }\; y\in Y(x),\;l=1,\ldots,L,\;\mbox{ and }\;p(x,u)-p_*(u)\le 0\big\}.
\end{equation}
Denote $f_{l+1}(x,y):=p(x,\bar u)-p_*(\bar u)$ and define
\begin{equation}\label{phi-sip}
\phi(x):=\max\Big\{p(x,\bar u)-p_*(\bar u),\;\max_{y\in Y(x),\atop 1\le l\le L}f_l(x,y)\Big\}.
\end{equation}
Applying the error bound from the second part of Theorem~\ref{th:local} to the system in \eqref{phi-sip0} and \eqref{phi-sip} with $L$ replaced by $L+1\ge 2$, we find $c_0>0$ such that
\begin{equation*}
{\rm dist}\big(x,S(\bar{u})\big)\le c_0\,\big[\phi(x)\big]_+^{\tau}\;\mbox{ for all }\;x\in G\;\mbox{ with }\;\tau={\mathscr R}\big(2n+(m+r+s+2)(n+1),d+L+2\big)^{-1}.
\end{equation*}
Since $p(x,\cdot)$ is locally Lipschitzian around $\ou$ and $G$ is compact, we find constants $\ell,\ve>0$ such that
\begin{equation}\label{eq:assumption2}
|p(x,u_1)-p(x,u_2)|\le\ell\|u_1-u_2\|\;\mbox{ for all }\;x\in G\;\mbox{ and }\;u_1,u_2\in {\mathbb B}_\ve(\bar u).
\end{equation}
Denote $c:=(2\beta^{-1}\ell)^{\tau}$ with $\beta:=c_0^{-\frac{1}{\tau}}> 0$ and for any $w\in S(u)$ select $z\in S(\bar {u})$ with $\|w-z\|={\rm dist}(w, S(\bar{u}))$. To justify the H\"older stability \eqref{Su}, it remains to verify the estimate
\begin{equation}\label{eq:claim}
\|w-z\|\le c\,\|u-\bar{u}\|^{\tau}.
\end{equation}
We consider the case of $w\in S(u)\backslash S(\bar{u})$; otherwise, (\ref{eq:claim}) is trivial. Then
\begin{eqnarray*}
p(w,\bar{u})-p_*(\bar{u})=\phi(w)=\big[\phi(w)\big]_+\ge\beta d\big(w,S(\bar{u})\big)^{\frac{1}{\tau}}=\beta \|w-z\|^{\frac{1}{\tau}}.
\end{eqnarray*}
Since $z\in S(\bar{u})$, we get $p(z,\bar{u})=p_*(\bar{u})\le p(w,\bar{u})$, and hence
\begin{equation}\label{eq:ok}
\|w-z\|^{\frac{1}{\tau}}\le\beta^{-1}\left(p(w,\bar{u})-p_*(\bar{u})\right)\le\beta^{-1}\big(p(w,\bar{u})-p(z,\bar{u})\big).
\end{equation}
Furthermore, it follows from $w\in S(u)$ that $p(w,u)\le p(z,u)$, and so (\ref{eq:assumption2}) gives us the relationships
\begin{eqnarray*}
p(w,\bar{u})-p(z,\bar{u})&=&\big(p(w,u)-p(z,u)\big)+\big(p(z,u)-p(z,\bar{u})\big)+\big(p(w,\bar{u})-p(w,u)\big)\\
&\le& \big(p(z,u)-p(z,\bar{u})\big)+\big(p(w,\bar{u})-p(w,u)\big)\le2\ell\|u-\bar{u}\|\;\mbox{ as }\;w,z\in G
\end{eqnarray*}
implying together with (\ref{eq:ok}) that $\|w-z\|^{\frac{1}{\tau}}\le\beta^{-1}\big(p(w,\bar{u})-p(z,\bar{u})\big)\le 2\beta^{-1}\ell \|u-\bar{u}\|$. Therefore
$$
{\rm dist}\big(w,S(\bar{u})\big)=\|w-z\|\le c\,\|u-\bar{u}\|^{\tau},
$$
which justifies (\ref{eq:claim}) and thus completes the proof of the theorem.\qed
\end{proof}\vspace*{-0.1in}
\begin{remark}{\bf (sharper exponents in higher-order stability for equality descriptions.)}\label{remark:2} Based on Remark~\ref{remark:21}, observe from the proof of Theorem~\ref{th:stability_SIP} that the H\"older stability exponent in \eqref{Su} can be sharpened to $\tau={\mathscr R}(2n+(m+r+s+2)(n+1),d+L+1)$ if the sets $Y(x)$ in $(GSIP)$ are described by the polynomial equalities only: $Y(x)=\{y|\;h_j(x,y)=0,\;j=1,\ldots,s\}$.\end{remark}\vspace*{-0.2in}

\subsection{Optimization of Polynomial Matrix Inequalities}\vspace*{-0.1in}

Next we address stability issues in optimization of {\em polynomial matrix inequalities}:
\begin{eqnarray*}
(PMI)\qquad\mbox{minimize }\;p_0(x)\;\mbox{ subject to }\;P(x)\preceq 0,
\end{eqnarray*}
where $p_0\colon\mathbb{R}^n\rightarrow\mathbb{R}$ is a real polynomial of degree not exceeded $d$, and where the mapping $P\colon\mathbb{R}^n\rightarrow S^m$ is such that each $(i,j)$th element of its image $\big(P(x)\big)_{ij}$, $1\le i,j\le m$, is a real polynomial under the same degree assumption. As in Subsection~5.1, consider the perturbed version of $(PMI)$:
\begin{eqnarray*}
(PMI_u)\qquad\mbox{minimize }\;p(x,u)\;\mbox{ subject to }\;P(x)\preceq 0,
\end{eqnarray*}
where $u\in\mathbb{R}^q$, $p(x,0)=p_0(x)$, the degree of $p(\cdot,u)$ is at most $d$, and $p(x,\cdot)$ is locally Lipschitzian.\vspace*{-0.1in}

\begin{theorem}{\bf (higher-order stability of optimal solution sets to polynomial matrix inequalities).}\label{pmi} Assuming that the feasible set $G:=\{x|\; P(x)\preceq 0\}$ in $(PMI_u)$ is nonempty and compact and then denoting by $S(u)$ the set of optimal solutions to this problem, for each $\ou\in\R^q$ we find constants $c,\ve>0$ such that
$$
S(u)\subset S(\bar u)+c\,\|u-\bar u\|^\tau\;\mbox{ if }\;\|u-\ou\|\le\ve\;\mbox{ with }\;\tau={\mathscr R}\big(2n+(m+3)(n+1),d+4\big)^{-1}.
$$
\end{theorem}\vspace*{-0.1in}
\begin{proof}
Consider the function $f(x,y):=\sum_{i,j=1}^m y_i\left(P(x)\right)y_j$, which is a polynomial of degrees at most $d+2$. We clearly have the equivalence
\[
P(x)\preceq 0\Longleftrightarrow f(x,y)\le 0\;\mbox{ for all }\;\|y\|^2=1.
\]
Defining now the constant sets $Y(x)$ by
$$
Y(x)\equiv\Omega:=\big\{y\in\mathbb{R}^m\big|\;\|y\|^2=1\big\},
$$
observe that the set $\O$ is regular. Then applying Theorem~\ref{th:stability_SIP} with $L=1$, $r=0$, and $s=1$, replacing $d$ with $d+2$, and taking into account Remark~\ref{remark:2} justify the claimed H\"older stability.\qed
\end{proof}\vspace*{-0.2in}

\subsection{Second-Order Cone Programs with Polynomial Data}\vspace*{-0.1in}

The class of polynomial optimization problems considers below belongs to {\em second-order cone programming}
\begin{eqnarray*}
(SOCP)\qquad\mbox{minimize }\;p_0(x)\;\mbox{ subject to }\;\big\|\big(f_1^{(l)}(x),\ldots,f_{m}^{(l)}(x)\big)\big\|\le f_{m+1}^{(l)}(x),\quad l=1,\ldots,L,
\end{eqnarray*}
where $p_0,f_j^{(l)}\colon\mathbb{R}^n\rightarrow\mathbb{R}$ for $j=1,\ldots,m+1$ and $l=1,\ldots,L$ are real polynomials with degrees at most $d$. Besides other areas of applications, problems of this type arise in optimization under {\em uncertainty} via the {\em robust optimization approach}. As a simple illustration, let us consider the following quadratic optimization problem under uncertainty in constraint data:
\begin{eqnarray}\label{robust}
\mbox{minimize }\;x^TAx+a^Tx\;\mbox{ subject to }\;x^TB_lx+b_l^Tx+\beta_l\le 0,\;l=1,\ldots,L,
\end{eqnarray}
where the triples $(B_l,b_l,\beta_l)\in S^n\times\mathbb{R}^n\times\mathbb{R}$, $l=1,\ldots,L$ are not known in advance. As commonly accepted, suppose that the constraint data belong to the {\em ellipsoidal uncertainty set}
\[
\mathcal{U}_l:=\Big\{(B_l^{(0)},b_l^{(0)},\beta_l^{(0)})+\sum_{j=1}^sv_l^j(B_l^{(j)},b_l^{(j)},\beta_l^{(j)})\Big|\;\|(v_l^1,\ldots,v_l^s)\|\le 1\Big\}.
\]
The robust optimization approach deals with an associated problem of minimizing an adverse effect of uncertainty. In our case it reads as follows:
\begin{eqnarray*}
\mbox{minimize }\;x^TAx+a^Tx\;\mbox{ subject to }\;x^TB_lx+b_l^Tx+\beta_l\le 0\;\mbox{ for all }\;(B_l,b_l,\beta_l)\in\mathcal{U}_l,\;l=1,\ldots,L.
\end{eqnarray*}
By direct calculations we have for each $l=1,\ldots,L$ that
\begin{eqnarray*}
&&\sup_\big{(B_l,b_l,\beta_l)\in\mathcal{U}_l}\{x^TB_lx+b_l^Tx+\beta_l\big\}\\
&=&x^TB_l^{(0)}x+(b_l^{(0)})^Tx+\beta_l^{(0)}+\sup_{\|(v_l^1,\ldots,v_l^s)\|\le 1}\sum_{j=1}^s v_l^j\big(x^TB_l^{(j)}x+(b_l^{(j)})^Tx+\beta_l^{(j)}\big)\\
&=&x^TB_l^{(0)}x+(b_l^{(0)})^Tx+\beta_l^{(0)}+\big\|\big(x^TB_l^{(1)}x+(b_l^{(1)})^Tx+\beta_l^{(1)},\ldots,x^TB_l^{(s)}x+(b_l^{(s)})^Tx+\beta_l^{(s)}\big)\big\|.
\end{eqnarray*}
This therefore reduces the uncertainly model to the above framework of (SOCP).\vspace*{-0.1in}

Now we consider the following perturbations of $(SOCP)$ given for all $u\in\R^q$ by
\begin{eqnarray*}
(SOCP_u)\qquad\mbox{minimize }\;p(x,u))\;\mbox{ subject to }\;\big\|\big(f_1^{(l)}(x),\ldots,f_{m}^{(l)}(x)\big)\big\|\le f_{m+1}^{(l)}(x),\quad l=1,\ldots,L,
\end{eqnarray*}
where $p(x,0)=p_0(x)$, $p(\cdot,u)$ is a polynomial of degree at most $d$, and $p(x,\cdot)$ is locally Lipschitzian on $\R^q$.\vspace*{-0.1in}

\begin{theorem}{\bf (higher-order stability of optimal solution sets to polynomial second-order cone programs).} Assuming that the feasible set
$G=\{x|\;\|(f_1^{(l)}(x),\ldots,f_{m}^{(l)}(x))\|\le f_{m+1}^{(l)}(x)\}$ to $(SOCP_u)$ is nonempty and compact, for any $u\in\R^q$ denote by $S(u)$ the set of optimal solutions to this problem. Then given $\ou\in\R^q$, there exist constants $c,\ve>0$ such that
\begin{equation}\label{s-socp}
S(u)\subset S(\bar u)+c\,\|u-\bar u\|^\tau\;\mbox{ if }\;\|u-\ou\|\le\ve\;\mbox{ with }\;\tau={\mathscr R}\big(2n+(mL+L+2)(n+1),d+L+1\big)^{-1}.
\end{equation}
\end{theorem}\vspace*{-0.1in}
\begin{proof}
For each $l=1,\ldots,L$ we define the function
\[
f_l(x,y):=\sum_{j=1}^{m}y_j^{(l)}f_j^{(l)}(x)-f_{m+1}^{(l)}(x),
\]
which is a polynomial with real coefficients on $\mathbb{R}^{n+m}$ of degree at most $d+1$. Note that
\[
\big\|\big(f_1^{(l)}(x),\ldots,f_{m}^{(l)}(x)\big)\big|\le f_{m+1}^{(l)}(x)\Longleftrightarrow\sup\Big\{\sum_{j=1}^{m}y_j^{(l)}f_j^{(l)}(x)-f_{m+1}^{l}(x)\Big|\; \big\|(y^{(l)}_1,\ldots,y^{(l)}_{m})\big\|=1\Big\}\le 0
\]
whenever $l=1,\ldots,L$. Consider the sets $Y(x)\subset\R^{mL}$ given by
$$
Y(x)=\Omega:=\big\{y=(y^{(1)},\ldots,y^{(L)})\big|\;y^{(l)}=(y^{(l)}_1,\ldots,y^{(l)}_{m}),\;\|y^{(l)}\|^2=1,\;l=1,\ldots,L\big\},\quad x\in\R^n.
$$
Observe that these sets are nonempty and uniformly bounded around any $x\in\R^n$ and that the set $\Omega$ is regular. Applying now Theorem~\ref{th:stability_SIP} and Remark~\ref{remark:2} with $r=0$ and $s=L$ and with replacing $m$ by $mL$, we arrive at the claimed H\"older stability in \eqref{s-socp}.\qed
\end{proof}\vspace*{-0.3in}

\section{Convergence Rates in Feasibility and Asymptotic Analysis}\vspace*{-0.15in}
\setcounter{equation}{0}

The concluding section of the paper presents applications of the obtained error bounds to deriving {\em explicit convergence rate estimates} in two important frameworks involving polynomial functions. The first issue concerns the cyclic projection algorithm to solve the well-known convex feasibility problem, while the second one deals with asymptotic analysis of subgradient dynamical systems.\vspace*{-0.2in}

\subsection{Cyclic Projection Method for Convex Feasibility Problems}\vspace*{-0.1in}

A common problem arising in diverse areas of mathematics and applications is to find a point in the intersection of closed convex sets $C_l$, $l=1,\ldots,L$. This problem is often referred to as the {\em convex feasibility problem}. One of the most powerful methods for solving the convex feasibility problem is the so-called {\em cyclic projection algorithm}, which is formulated as follows. Given finitely many closed convex sets $C_1,\ldots C_L$ in $\mathbb{R}^n$ with $\bigcap_{l=1}^L C_l\ne\emp$, pick $x_0\in\mathbb{R}^n$ and denote $\pi_l:=\pi_{C_l}$ for $l=1,\cdots,L$, where $\pi_{C_l}$ stands for the (unique) Euclidean projection to the set $C_l$. The sequence of \emph{cyclic projections} $(x_k)_{k\in\NN}$ is defined by
\begin{equation}\label{Cycse:1}
x_1:=\pi_1 x_0,\;\ldots,x_L:=\pi_Lx_{L-1},\;x_{L+1}:=\pi_1x_{L},\ldots.
\end{equation}
In the case of $L=2$ the cyclic projection method reduces to the well-known {\em von Neumann alternating projection method}; see, e.g., \cite{Heinz,BC2011} and the references therein, where the reader can find various practical applications of the cyclic projection method and related algorithms.\vspace*{-0.1in}

Convergence properties of the cyclic projection algorithm have been examined by many researchers. In particular, the seminal paper by Bregman \cite{Bregman} shows that the sequence $(x_k)_{k\in\NN}$ generated by cyclic projections always converges to a point in $C$. Furthermore, it is proved in \cite{GPR} that the cyclic projection algorithm converges linearly for regular intersections of convex sets. On the other hand, for convex sets with irregular intersections (e.g., when the intersection is a singleton), the cyclic projection algorithm may not exhibit linear convergence even in simple two-dimensional cases as observed by \cite[Example~5.3]{BBJAT1}. This raises the following basic question: Can we estimate the convergence rate of the cyclic projection algorithm for convex sets with possibly {\em irregular} intersections of certain particular structures? Quite recently \cite{Borwein} this issue has been investigated in the setting where the convex sets $C_l$ are the so-called {\em basic semialgebraic} given by
\begin{equation}\label{basic-sa}
C_l:=\big\{x\in\mathbb{R}^n\big|\;g_i^l(x)\le 0,\;i=1,\ldots,r_l\big\},\quad l=1,\ldots,L,
\end{equation}
with convex polynomials $g_i^l$ of degrees at most $d$. It is shown in \cite{Borwein} that for $d>1$ the sequence of cyclic projections $(x_k)_{k\in\NN}$ from \eqref{Cycse:1} converges at the rate of
$$
\disp\frac{1}{k^{\rho}}\;\mbox{ with }\;\rho:=\disp\frac{1}{\min\left\{(2d-1)^n-1,\,2{n-1\choose{[(n-1)/2]}}d^n-2\right\}},
$$
where the symbol $[a]$ stands for the integer part of the number $a$\vspace*{-0.1in}.

In what follows we establish an explicit convergence rate estimate for the more general case of convex sets described by polynomial matrix inequalities. It is said that a set $C$ is {\em polynomial matrix inequality representable} if there are numbers $m,d\in\mathbb{N}$ and a mapping $P\colon\mathbb{R}^n\rightarrow S^m$ such that each $(i,j)$th element of its image $\big(P(x)\big)_{ij}$, $1\le i,j\le m$ is a real polynomial with degree at most $d$ and we have the representation
$C=\{x\in\mathbb{R}^n|\;P(x)\preceq 0\}$. It is obvious that any basic semialgebraic set \eqref{basic-sa} belongs to this class. The following example shows that
the converse is not true.\vspace*{-0.1in}

\begin{example} {\bf (polynomial matrix inequality representable sets may not be basic semialgebraic).}\label{not-basic}
Given an odd number $d$, consider the set
$$
C:=\left\{(x_1,x_2)\in\mathbb{R}^2\Big|\;\left(\begin{array}{cc}
-x_1^d&1\\
1&-x_2^d
\end{array}\right)\preceq 0\right\}=\big\{(x_1,x_2)\big|\;x_1x_2\ge 1,\;x_1\ge 0,x_2\ge 0\big\},
$$
which is clearly convex and polynomial matrix inequality representable. Let us show that $C$ cannot be written in the basic semialgebraic form \eqref{basic-sa}. Assuming the contrary, take the function $f(x_1,x_2):=x_1$ for which $\inf_{x=(x_1,x_2)\in C}f(x)=0$. It follows from the polynomial extension of the celebrated (quadratic) Frank-Wolfe theorem that any convex polynomial $f$ bounded from below over a basic semialgebraic convex set attains its minimum; see \cite{convex_polynomial}. This gives us $\bar x\in C$ such that $f(\bar x)=\inf_{x\in C}f(x)=0$ for the function $f$ and the set $C$ defined above. It surely contradicts the structures of $f$ and $C$, and thus verifies the claim.
\end{example}\vspace*{-0.15in}

To derive an explicit convergence rate estimate for the cyclic projection algorithm applied to polynomial matrix inequality representable convex sets, we use the error bound result from Theorem~\ref{th:local} and the following abstract result on cyclic projections established in \cite{Borwein} under a certain H\"older condition on convex sets.\vspace*{-0.1in}

\begin{lemma}{\bf (conditional cyclic convergence rate).}\label{TheMTR:1} Let $C_l\subset\R^n$, $l=1,\ldots,L$, be closed convex sets with $C:=\bigcap_{l=1}^L C_l\ne\emp$, and let $(x_k)_{k\in\NN}$ be the sequence of cyclic projections \eqref{Cycse:1} with a starting point $x_0$. Impose the following H\"older condition of rank $\tau\in(0,1]$: for any compact set $K\subset\mathbb{R}^n$ there is a positive number $c_0$ providing the estimate
\begin{equation}\label{remark:3}
di(x,C)\le c_0\bigg(\sum_{l=1}^L \di(x,C_l)\bigg)^{\tau}\;\mbox{ whenever }\;x\in K.
\end{equation}
Then there exists $x_{\infty}\in C$ such that the sequence $(x_k)_{k\in\NN}$ converges to $x_\infty$, and for some $M>0$ we have
\begin{eqnarray*}
\|x_{k}-x_{\infty}\|\le M\,k^{-\dfrac{1}{2\tau^{-1}-2}}.
\end{eqnarray*}
\end{lemma}\vspace*{-0.1in}

Now we are now ready to obtain an explicit convergence rate estimate of the cyclic projection algorithm for polynomial matrix inequality representable convex sets.\vspace*{-0.1in}

\begin{theorem}{\bf (explicit convergence rate of the cyclic projection algorithm).}\label{TheMTR:3} Let $C_l$ be polynomial matrix inequality representable convex sets given in the form
\[
C_l:=\big\{x\in\mathbb{R}^n\big|\;A^{(l)}(x)\preceq 0\big\},\quad l=1,\ldots,L,
\]
where every $A^{(l)}\colon\mathbb{R}^n\rightarrow S^m$ is a matrix mapping such that each $\big(A^{(l)}(x)\big)_{ij}$, $1\le i,j\le m$, is a real polynomial with degree at most $d\ge 1$. Given $x_0\in\R^n$, consider the sequence $(x_k)_{k\in\NN}$ of cyclic projections {\rm (\ref{Cycse:1})}. Then $x_k$ converges to some $x_{\infty}\in C$ and there is $M>0$ ensuring the estimate
\begin{align*}
\|x_{k}-x_{\infty}\|\le M\,\frac{1}{k^{\rho}},\quad k\in\mathbb{N},
\end{align*}
where $\rho:=[2\,{\mathscr R}(2n+(m+3)(n+1),{d}+L)-2]^{-1}$ with the constant $\mathscr R$ defined in \eqref{loj}.
\end{theorem}\vspace*{-0.1in}
\begin{proof}
Let us first show that for any compact set $K\subset\mathbb{R}^n$ there is $c_0>0$ such that condition \eqref{remark:3} holds with some explicitly calculated $\tau\in(0,1)$. To proceed, fix a compact set $K$ and define
$$
\Omega:=\big\{y\in\mathbb{R}^m\big|\;\|y\|^2=1\big\}\;\mbox{ and }\;\disp f_l(x,y):=\sum_{i,j=1}^m y_i\left(A^{(l)}(x)\right)y_j,\;l=1,\ldots,L.
$$
Then we readily have the polynomial set representations
\begin{equation}\label{poly-C}
C_l=\big\{x\in\mathbb{R}^n\big|\;f_l(x,y)\le 0\;\mbox{ for all }\;y \in\Omega\big\}\;\mbox{ and }\;C=\Big\{x\in\mathbb{R}^n\Big|\;\disp\max_{y\in\Omega\atop 1\le l\le L}f_l(x,y)\le 0\Big\}.
\end{equation}
Applying now the error bound result of Theorem~\ref{th:local} to the set $C$ from \eqref{poly-C} with $Y(x)\equiv\Omega$, $r=0$, and $s=1$ while taking into account Remark~\ref{remark:21} and the regularity of the set $\O$ above, we find $c_0$ such that
\begin{eqnarray*}
{\rm dist}(x,C)\le c_0\,\bigg(\big[\max_{y\in\Omega\atop 1\le l\le L}f_l(x,y)\big]_{+}\bigg)^{\tau}\le c_0\,\bigg(\sum_{l=1}^L\big[\max_{y\in\Omega}f_l(x,y)\big]_{+}\bigg)^{\tau}\;\mbox{ for all }\;x\in K,
\end{eqnarray*}
where $\tau={\mathscr R}(2n+(m+3)(n+1),{d}+L)^{-1}$.
For each $l=1,\ldots,L$ define the function $h_l(x):=[\max_{y\in\Omega}f_l(x,y)]_+$, $x\in\mathbb{R}^n$, and observe that it is locally Lipschitzian by Proposition~\ref{Lemma31}. Due to the compactness of $K$ we have $M_0:=\max_{x\in K}\sum_{l=1}^L d(x,C_l)<\infty$. Since the set $\hat{K}:=K+M_0\overline{\mathbb{B}}$ is also compact, we get $c_l>0$ ensuring the global Lipschitz condition on this set:
\[
|h_l(x_1)-h_l(x_2)|\le c_l\|x_1-x_2\|\;\mbox{ for all }\;x_1,x_2\in\hat{K},\quad l=1,\ldots,L.
\]
Now pick an arbitrary vector $x\in K\subset\hat K$ and select $\hat{x}_l\in C_l$ such that ${\rm dist}(x,C_l)=\|x-\hat{x}_l\|$. Denoting $c:=\max\{c_1,\ldots,c_L\}>0$ and observing that $\hat{x}_l\in\hat{K}$ and $h_l(\hat{x}_l)=0$ give us the relationships
\[
\sum_{l=1}^L{\rm dist}(x,C_l)=\sum_{l=1}^L\|x-\hat{x}_l\|\ge\sum_{l=1}^L\frac{1}{c_l}|h_l(x)-h_l(\hat{x}_l)|\ge\frac{1}{c}\sum_{l=1}^L|h_l(x)|=\frac{1}{c}\sum_{l=1}^L h_l(x),
\]
which ensure therefore the estimates
\[
{\rm dist}(x,C)\le c_0\,\bigg(\sum_{l=1}^L h(x)\bigg)^{\tau}\le c_0\,\bigg(c\,\sum_{l=1}^L{\rm dist}(x,C_l)\bigg)^{\tau}.
\]
Thus we arrive at (\ref{remark:3}) and complete the proof of the theorem by using Lemma~\ref{TheMTR:1}.\qed
\end{proof}\vspace*{-0.2in}

\subsection{Subgradient Dynamical Systems}\vspace*{-0.1in}

The final piece of this paper is devoted to applications of the obtained error bounds to conducting a quantitative asymptotic analysis of the {\em subgradient dynamical systems} given by
\begin{eqnarray*}
(DS)\qquad\begin{array}{cc}
&\quad\;\;0\in\dot{x}(t)+\partial\phi\big(x(t)\big)\;\mbox{ a.e. on }\;(0,T),\\
&\partial\phi\big(x(t)\big)\ne\emp\;\mbox{ for all }\;t\in[0,T),
\end{array}
\end{eqnarray*}
where $\phi$ is a semialgebraic function on $\mathbb{R}^n$ in the form $\phi(x):=\sup\{f(x,y)|\;y\in\Omega\}$. We suppose in what follows that $\Omega$ is a regular, compact, and semialgebraic set defined by
\begin{equation}\label{eq:Omegaf}
\Omega:=\big\{y\in\mathbb{R}^m\big|\;g_i(y)\le 0\;\mbox{ for }\;i=1,\ldots,r\;\mbox{ for }\;h_j(y)=0,\;j=1,\ldots,s\}
\end{equation}
and that all the functions $f$, $g_i$, and $h_j$ are real polynomials with degrees at most $d$. By a solution of $(DS)$ we understand any absolutely
continuous curve $x\colon[0,T)\rightarrow\mathbb{R}^n$ satisfying the relationships therein. A trajectory $x(t)$ is called {\em maximal} if there is no possible extensions of its domain compatible with $(DS)$. Our asymptotic analysis of trajectories for $(DS)$ as $t\to\infty$ is based on Theorem~\ref{LojasiewiczGradientTheorem} extending the \L ojasiewicz inequality to maximum functions and on the following result taken from \cite{Lewis}.\vspace*{-0.1in}

\begin{lemma} {\bf (convergence to critical points).}\label{crit} Let $x(t)$ be a bounded maximal trajectory of $(DS)$, where the function $\phi$ is assumed to be bounded from below on $\mathbb{R}^n$. Then  $x(t)$ is defined on $[0,\infty)$ and converges to a critical point $\ox\in\mathbb{R}^n$ of $\phi$, i.e., such that $0\in\partial\phi(\ox)$. Take any $\theta\in[0,1)$ for which the quantity $\frac{|\phi(x)-\phi(\ox)|^{\theta}}{m_f(x)}$ is bounded for any $x$ around $\ox$ under the convention that $\frac{0}{0}=1$. Then there are numbers $c,t_0>0$ such that
\[
\|x(t)-\ox\|\le c\,(t+1)^{-(\frac{1-\theta}{2\theta-1})}\;\mbox{ for all }\;t\ge t_0.
\]
\end{lemma}\vspace*{-0.1in}

We finally arrive at the constructive rate evaluation of asymptotic convergence for maximal trajectories.\vspace*{-0.1in}

\begin{theorem}{\bf (explicit rate estimate of asymptotic convergence).}\label{asym} Let $\phi$ be a semialgebraic function on $\mathbb{R}^n$ given by $\phi(x):=\max\{f(x,y)|\;y\in\Omega\}$, where $\Omega$ is a semialgebraic set in form \eqref{eq:Omegaf} under the assumptions formulated above. Suppose in addition that $\phi$ is bounded from below on $\mathbb{R}^n$. Then any bounded maximal trajectory $x(t)$ for $(DS)$ converges as $t\to\infty$ to a critical point $\ox\in\mathbb{R}^n$ of $\phi$. Furthermore, there are numbers $c,t_0>0$, which ensure the estimates
\begin{equation}\label{DS-rate}
\|x(t)-\ox\|\le c\,(t+1)^{-(\frac{1-\theta}{2\theta-1})}\;\mbox{ and }\;|\phi(x(t))-\phi(\ox)|\le c\,t^{-\frac{1}{2\theta-1}}\mbox{ for all }\;t \ge t_0,
\end{equation}
where $\theta=1-{\mathscr R}(2n+(m+r+s)(n+1),d+2)^{-1}$ with the constant $\mathscr R$ defined in \eqref{loj}.
\end{theorem}\vspace*{-0.1in}
\begin{proof}
Let $\ox$ be a critical point of $\phi$ the existence of which follows from Lemma~\ref{crit}. Denote $\tau:={\mathscr R}(2n+(m+r+s)(n+1),d+2)^{-1}$ and $\theta:=1-\tau$. Then Theorem~\ref{LojasiewiczGradientTheorem} tells us that there are numbers $c_0,\ve>0$ such that ${\frak m}_{\phi}(x)\ge c|\phi(x)-\phi(\ox)|^{\theta}$ for all $x\in\overline{\mathbb{B}}_{\ve}(\ox)$, and therefore
\begin{equation}\label{eq:useful09}
\frac{|\phi(x)-\phi(\ox)|^{\theta}}{m_{\phi}(x)}\le c_0^{-1}\;\mbox{ whenever }\;x\in\overline{\mathbb{B}}_{\ve}(\ox).
\end{equation}
Thus the first estimate in \eqref{DS-rate} follows directly from Lemma~\ref{crit}. Define now $\gg\colon[0,\infty)\rightarrow\mathbb{R}$ by
$$
\gg(t):=\phi\big(x(t)\big)-\phi(\ox),\quad t\in[0,\infty),
$$
and observe that both functions $\phi$ and $\gg$ are Lipschitz continuous and so a.e.\ differentiable on $(0,\infty)$. This yields $\partial\phi(x(t))=\{\nabla\phi(t)\}$ for a.e.\ $t\in(0,\infty)$; see, e.g., \cite[Corollary~1.82]{Boris_Book}. Then we see that
\begin{equation}\label{eq:claim00}
\dot{\gg}(t)=\big\langle\nabla\phi(x(t)),\dot{x}(t)\big\rangle=-\big\|\nabla\phi(x(t))\big\|^2=-{\frak m}_{\phi}(x(t))^2\le-c_0^2 \, \big|\phi(x(t))-\phi(\ox)\big|^{2\theta}=-c_0^2\,\gg(t)^{2\theta}\le 0,
\end{equation}
where the second equality follows from $0\in\dot{x}(t)+\partial\phi(x(t))$ while the first inequality follows from \eqref{eq:useful09}. This shows, in particular, that the function $\gg(t)$ is nonincreasing. Since $x(t)\rightarrow\ox$ and $\gg(t)\to 0$ as $t\rightarrow\infty$, we get that either $\gg(t)\equiv 0$, or $\gg(t)>0$ for all large $t>0$. Integrating \eqref{eq:claim00} allows us to find $c,t_0>0$ such that the second estimate in \eqref{DS-rate} holds for all $t\ge t_0$, and thus we complete the proof of the theorem.\qed
\end{proof}
\end{document}